\documentclass[12pt]{amsart}
\usepackage{amsfonts}
\usepackage{amsmath}
\usepackage{amsxtra}
\usepackage{amssymb,latexsym}
\usepackage[mathcal]{eucal}
\usepackage{amscd}
\usepackage{graphicx}

\graphicspath{ {pdimages/} }

\newtheorem{theo}{{\bfseries Theorem}}[section]
\newtheorem{prop}[theo]{{\bfseries Proposition}}
\newtheorem{lem}[theo]{{\bfseries Lemma}}
\newtheorem{cor}[theo]{{\bfseries Corollary}}
\newtheorem{df}[theo]{{\bfseries Definition}}

\newtheorem{ex}[theo]{{\bfseries Example}}

\def \ul {\underline}

\def \R {\mathbb R}

\def \I {\mathcal I}
\def \J {\mathcal J}

\def \S {\mathcal S}

\def \e {\epsilon}
\def \d {\delta}
\def \g {\gamma}

\def \ee {{\mathbf e}}
\def \pp {{\mathbf p}}
\def \qq {{\mathbf q}}
\def \vv {{\mathbf v}}

\def \MM {{\mathbf M}}

\def \SX {{\mathbf S_X}}
\def \SY {{\mathbf S_Y}}

\usepackage{amssymb,latexsym}
\usepackage[mathcal]{eucal}
\usepackage{amscd}
\usepackage{amsmath}
\usepackage{graphicx}
\usepackage{amsfonts}
\usepackage{amscd}

\numberwithin{equation}{section}

\renewcommand{\arraystretch}{1.5}

\begin{document}

\title[Smale vs. Markov]{  Good Strategies for the \\ Iterated Prisoner's Dilemma :\\ Smale vs. Markov}
\vspace{1cm}
\author{Ethan Akin}
\address{ Mathematics Department \\
    The City College \\
    137 Street and Convent Avenue \\
    New York City, NY 10031, USA }
\email{ethanakin@earthlink.net}

    \vspace{.5cm}
\date{March, 2017}

\begin{abstract}
In 1980 Steven Smale introduced a class of strategies for the Iterated Prisoner's Dilemma which used
as data the running average of the previous payoff pairs. This approach is quite different from the Markov
chain approach, common before and since, which used as data the outcome of the just previous play,
the \emph{memory-one} strategies. Our purpose here is to compare these two approaches focusing upon
\emph{good strategies} which, when used by a player, assure that the only way an opponent can obtain at
least the cooperative payoff is to behave so that both players receive the cooperative payoff. In addition,
we prove a version for the Smale approach of the so-called Folk Theorem concerning the existence of Nash equilibria
in repeated play. We also
consider the dynamics when certain simple Smale strategies are played against one another.
\end{abstract}



 \maketitle

 \tableofcontents

 \section{Introduction}
 \vspace{.5cm}

 The Iterated Prisoner's Dilemma has been an object of considerable study ever since Axelrod's description of the results
 of computer tournaments \cite{Ax} and Maynard Smith's application of game theory to evolutionary competition \cite{M-S}.
 Most of this work has focused upon what I call here \emph{Markov} strategies, i.e. memory-one plans. Competition between
 two such stategies
 leads to a Markov process on the set of outcomes, e.g.  \cite{BNS}, \cite{PD} and   various surveys cited below. Considerable
 simulation work has been done to analyze numerically the competition between such strategies, e.g.\cite{KYC} and \cite{SP}.
In this context, I characterized in \cite{A-13} and \cite{A-16} certain so-called \emph{good} strategies which ensured that
an opponent could receive at least the cooperative payoff only by playing so that both players receive exactly the cooperative
payoff. Recently, Mike Shub pointed out to me that Steve Smale had described in 1980 strategies which he called good
\cite{Sm}. These strategies use an entirely different way of aggregating the data of past outcomes. While there has been
some work on Smale's procedure, e.g. \cite{Abh}, \cite{BH} and \cite{BBH}, most of the game theory literature has ignored it.
For example, Smale's work is not referred to in \cite{KYC}, \cite{N}, \cite{HS}, \cite{S-93} or \cite{S-10}. Our purpose here
is to compare the Smale and Markov procedures and especially to use these to clarify the notion of a good strategy.
In addition, we use the Taylor-Jonker equations for evolutionary dynamics \cite{TJ} to analyze competition among
certain \emph{simple Smale strategies}.

The author would like to express his appreciation to the referee for his helpful comments and corrections.
\vspace{1cm}

 \section{Plans for the Iterated Prisoner's Dilemma}
 \vspace{.5cm}

 We will focus mostly on the symmetric version of the \emph{Prisoner's Dilemma}.
Each of the two players, X and Y, has a choice between
two strategies, $c$ and $d$. Thus, there are four outcomes which we list in the order: $cc, cd, dc, dd,$ where,
for example,  $cd$ is the outcome when X plays $c$ and Y plays $d$. Either player can use a \emph{mixed  strategy},
  randomizing by choosing $c$ with probability $p_c$ and
$d$ with the complementary probability $1-p_c$.

Each then receives a payoff. The following
 $2 \times 2$ chart describes the payoff to the X player. The transpose is the Y payoff.
\renewcommand{\arraystretch}{1.5}
\begin{equation}\label{1}
\begin{array}{|c||c|c|}\hline
X\backslash Y & \quad c \quad & \quad  d \quad \\ \hline \hline
c & \quad R \quad & \quad S \quad  \\ \hline
d & \quad T \quad  & \quad P \quad \\ \hline
\end{array}
\end{equation}

Alternatively, we can define the \emph{payoff vectors} for each player by
\begin{equation}\label{2}
\SX \quad = \quad (R, S, T, P) \qquad \mbox{and} \qquad \SY \quad = \quad (R, T, S, P).
\end{equation}

The payoffs are assumed to satisfy
\begin{equation}\label{3}
T \ > \ R  \ > \ P  \ > \ S \qquad \mbox{and} \qquad 2R  \ > \ T + S,
\end{equation}
but $2P$ might lie on either side of $T + S$.

In the Prisoner's Dilemma, the strategy $c$ is \emph{cooperation}.  When both players cooperate
they each receive the reward for cooperation (= $R$).
The strategy $d$ is \emph{defection}.  When both players defect they each receive the punishment for defection (= $P$).
However, if one player cooperates and the other does not, then the defector receives the large temptation payoff (= $T$), while
the hapless cooperator receives the very small sucker's payoff (= $S$). The condition $2R > T + S$ says that the reward
for cooperation is larger than the players would receive by dividing equally the total payoff of a $cd$ or $dc$ outcome.
Thus, the maximum total payoff occurs uniquely at $cc$ and that location is a \emph{strict Pareto optimum}, which means that
at every other outcome at least one player does worse. The cooperative
outcome $cc$ is clearly where the players ``should" end up. If they
could negotiate a binding agreement in advance of play, they would agree to play $c$ and each receive R. However, the structure
of the game is such that, at the time of play, each chooses a strategy in ignorance of the other's choice.
Furthermore, the strategy $d$ \emph{strictly dominates} strategy $c$.  This means that,
whatever Y's choice is, X receives a larger payoff by playing $d$ than by using $c$. In the array (\ref{1}) each number in the
$d$ row is larger than the corresponding number in the $c$ row above it. Hence, X chooses $d$, and for exactly the same reason,
Y chooses $d$. So they are driven to the $dd$ outcome with payoff P for each.

In the search for a way to avoid the mutually inferior payoff $(P,P)$,
attention has focused upon \emph{repeated play}.  X and Y play repeated rounds of the same game. For each round
the players' choices are made independently, but each is aware of all of the previous outcomes. The hope is that the threat of
future retaliation will rein in the temptation to defect in the current round. It is this \emph{Iterated Prisoner's Dilemma}
which we will consider here.

After the $k^{th}$ round the players receive payoffs $S^k = (S^k_X,S^k_Y)$ determined by the payoff matrix. For a single payoff
after $N$ rounds of play we use the time average:
\begin{align}\label{4}
s^N \quad = \quad  (s^N_X,s^N_Y) \quad = \quad \frac{1}{N} \ \Sigma_{k=1}^{N} \ S^k.
\end{align}

Observe that
\begin{equation}\label{4a}
s^{N+1} \quad = \quad \frac{N}{N+1} s^N \ + \ \frac{1}{N+1} S^{N+1}.
\end{equation}
and so
\begin{equation}\label{4b}
s^{N+1} - s^N \quad = \quad  \frac{1}{N+1} (S^{N+1} - s^N).
\end{equation}

The vector $s^N = (s^N_X,s^N_Y)$ lies in $\S$ defined to be the convex hull of the four payoff pairs.  If
$P < (T + S)/2$, then $\S$ is a quadrilateral with vertices: $(S,T), (R,R), (P,P), (T,S)$.  If
$P \geq  (T + S)/2$, then $\S$ is the triangle with vertices: $(S,T), (R,R), (T,S)$, which contains $(P,P)$.

The Euclidean
diameter of the set $\S$ is $\sqrt{2}(T-S)$. For $\e > 0$ we will call $N^*$ an $\e$ \emph{time-step} when $N^* > \sqrt{2}(T-S)/\e$.
In that case, for $N \geq N^*$ the distance $||s^{N+1} - s^N||$ of the step from $s^N$ to $s^{N+1}$ is less than $\e$.

\begin{prop}\label{prop00} For any infinite sequence of outcomes, the set $\Omega$ of limit points of the sequence
$\{ s^N \}$ of payoff pairs is closed, connected subset of $\S$. If $U$ is an open subset of $\S$ which contains $\Omega$, then
there exists $N^*$ such that $s^N \in U$ for all $N \geq N^*$.\end{prop}

{\bf Proof:} From (\ref{4b}) it is clear that $|| s^{N+1} - s^N || \to 0$ as $N \to \infty$. So the conclusion is
immediate from the following well-known lemma.

 $\Box$ \vspace{.5cm}

For the lemma we introduce a bit of notation.

For $A$ a nonempty, closed subset of a compact metric space $X$ and $x \in X$ we let $d(x,A) = \inf \{ d(x,a) : a \in A \}$.
If $B$ is another nonempty, closed subset we let $d(A,B) = \inf \{ d(a,b) : a \in A, b \in B \}$. By compactness,
there exists $a \in A$ such that $d(x,A) = d(x,a)$, a point in $A$ closest to $x$,  and
there exist $a \in A, b \in B$ such that $d(A,B) = d(a,b)$. Given $\e > 0$ we let $V_{\e}(A)$ denote the
open set $\{ x : d(x,A) < \e \}$, i.e.
the set of points less than $\e$ away from a point of $A$.

 \begin{lem}\label{lem00a} Let $\{ x^N \}$ be a sequence in a compact metric space $X$.  If $d(x^{N+1},x^N) \to 0$ as
 $N \to \infty$ then the set of limit points $\Omega$ is a nonempty, closed, connected subset of $X$.
 If $U$ is an open set containing $\Omega$ then
there exists $N^*$ such that $x^N \in U$ for all $N \geq N^*$.\end{lem}

 {\bf Proof:} The set limit points $\Omega$ is the intersection of the decreasing sequence $\{ X^k \}$ with $X^k$ the
 closure of the  tail ${\{ x^N : N \geq k \}}$ of the sequence. If $U$ is an open set containing $\Omega$ then
 $\{ U \}$ and $\{ X \setminus X^k : k = 1, \dots \}$ is an open cover of $X$ and so has a finite
 subcover. Since the $X^k$'s are decreasing, it follows that for some $N^*$, $\{U, X \setminus X^{N*} \}$ is a
 cover of $X$. Hence, $\{ x^N : N \geq N^* \} \subset U $. If $\Omega$ were empty we could apply this to $U = \emptyset$
 and get a contradiction.

 Now let $A_0$ and $A_1$ be disjoint nonempty, closed subsets of $\Omega$.
 We will see that $\Omega \setminus (A_0 \cup A_1)$ is nonempty and this implies that $\Omega$ is connected. Let $3 \e$ be the
 distance between the sets $d(A_0,A_1)$. Choose $n^*$ so that
 $d(x^{N+1},x^N) < \e$ for all $N \geq n^*$.

 The sequence repeatedly approaches arbitrarily closely to each point of $\Omega$.  Since $A_0$ and $A_1$ are nonempty we
 can define $n_0$ to be the minimum $n \geq n^*$ such that $x^n$ lies in  $V_{\e}(A_0)$.
 Inductively, for $k \geq 0$,
 define $n_{2k + 1}$ to be the minimum $n \geq n_{2k}$ such that $x_n  \in V_{\e}(A_1)$, and for $k \geq 1$, $n_{2k}$
  define $n_{2k}$ to be the minimum $n \geq n_{2k - 1}$ such that $x_n  \in V_{\e}(A_0)$. It is clear
  that $x^{n_i - 1} \in V_{2 \e}(A_0) \setminus V_{\e }(A_0)$ if $i$ is even and is in
  $V_{2 \e}(A_1) \setminus V_{\e }(A_1)$ if $i$ is odd. Hence, the subsequence $\{  x^{n_i - 1} \}$ lies
  in the closed set $X \setminus (V_{\e }(A_0) \cup V_{\e }(A_1))$ and so it has limit points not in $A_0 \cup A_1$.

 $\Box$ \vspace{.5cm}

The choice of play for the first round is the \emph{initial play}.
A \emph{strategy}  is a choice of initial play together with what we will call a \emph{plan}.
A plan is a choice of play, after the first round, to respond to any possible past history of outcomes
in the previous rounds.

If X and Y use strategies with only pure strategy choices then the result is an infinite sequence
of outcomes. However, if mixed strategy choices are used either on the initial plays or as part of the plan
there results a probability measure on the space of all such sequences.

We will consider plans which use just a crucial portion of the past history data.

We will use the label \emph{Markov plan} for  a \emph{ stationary, memory-one plan}
 which bases its response entirely on the outcome of the previous round.
 For example, the \emph{Tit-for-Tat} plan, hereafter $TFT$, due to Anatol Rapoport,
 plays the opponent's response from the previous play.

With the outcomes listed in order as $cc, cd, dc, dd$, a Markov plan for X is given by a vector
$\pp = (p_1, p_2, p_3, p_4) = (p_{cc},p_{cd},p_{dc},p_{dd})$ where $p_z$
is the probability of playing c when the  outcome $z$ occurred in the previous round.
On the other hand, if Y also uses a Markov plan
$\qq = (q_1,q_2,q_3,q_4) $ then
the  response vector is $(q_{cc},q_{cd},q_{dc},q_{dd}) = (q_1, q_3, q_2, q_4)$ and
 the successive outcomes follow a Markov chain with transition matrix
given by:
\begin{equation}\label{5}
\MM \quad =\quad \begin{pmatrix}p_1q_1 & p_1(1-q_1) & (1-p_1)q_1 & (1-p_1)(1-q_1)  \\
p_2q_3 & p_2(1-q_3) & (1-p_2)q_3 & (1-p_2)(1-q_3) \\
p_3q_2 & p_3(1-q_2) & (1-p_3)q_2 & (1-p_3)(1-q_2) \\
p_4q_4 & p_4(1-q_4) & (1-p_4)q_4 & (1-p_4)(1-q_4)
\end{pmatrix}.
\end{equation}
Notice the switch in numbering from the Y strategy $\qq $ to the Y response vector. This is done because switching the perspective
of the players interchanges $cd$ and $dc$. This way
the ``same" plan for X and for Y is given by the same vector.  For example, $TFT$ for X and for Y
is given by $\pp  = \qq = (1,0,1,0)$, but the response vector for Y is $ (1,1,0,0)$.  The plan   \emph{Repeat}
is given by $\pp   = \qq = (1,1,0,0)$ with the response vector for Y  equal to $ (1,0,1,0)$. This plan just repeats
the previous play, regardless of what the opponent did.

We can think of a Markov chain on a finite set $I$ (in this case $I = \{ cc, cd, dc, dd \})$ as representing motion
on a directed graph with vertices $I$ and an edge from $i_1$ to $i_2$ if there is a positive probability, according to
$\MM$, of moving from $i_1$ to $i_2$, i.e.  $\MM_{i_1 i_2} > 0$. In particular, there is an edge from $i$ to itself
 when $\MM_{i i} > 0$.  A \emph{path} in the graph is a state sequence $i^1,...,i^n$ with
 $n > 1$ such that there is an edge from $i^k$ to $i^{k+1}$ for $k = 1,...,n-1$.  A set of
 states $J \subset I$ is called a \emph{closed set} when it is  nonempty  and
no path  that begins in $J$ can
exit $J$.  For example, the entire set of states is closed and for any  $i$ the
set of states accessible via a path that begins at $i$ is a closed set. The subset $J $ is called a \emph{terminal set}
when it is closed and when for $i_1, i_2 \in J$ there exists a path from $i_1$ to $i_2$.  Equivalently, a
terminal set is a minimal closed set. Since the set is closed,
the path moves only on elements of $J$.

A vector $\vv$ is a \emph{stationary distribution} for $\MM$ when $\vv_i \geq 0, \Sigma_i \vv_i = 1$
and $\vv \MM = \vv$. For a terminal set $J$ there is a unique stationary distribution vector $\vv_J$ such that
$(\vv_J)_i = 0 $ for $ i \not\in J$.
Furthermore,
$(\vv_J)_i  > 0 $ for $ i \in J$.
Restricted to a terminal set, the system is ergodic and so for any function $f : I \to \R$ the sequence of
time averages $\{ \frac{1}{N} \Sigma_{k = 1}^{N} f(i^k) : T = 1, 2, \dots \}$ converges with probability $1$
to the space average $\Sigma_{i \in J}  (\vv_J)_i f(i)$, the expected value with respect to $\vv_J$.
 That is, such convergence occurs except on a set of outcome sequences which has probability $0$.
Think of the outcome sequence for a fair coin such that Heads comes up every time.

Distinct terminal sets are disjoint. If $i \in I$ lies in some terminal set then it is called \emph{recurrent};
otherwise, it is called \emph{transient}. If $i$ is transient then for each terminal set $J$ there is a probability,
determined by $i$ and $\MM$ so that the process beginning from $i$ eventually enters $J$. This probability might be zero
but when summed over all of terminal sets the probabilities add up to $1$. If the process enters $J$ then the
time average for any function $f$ approaches the expected value with respect to $\vv_J$ with probability $1$.

The matrix $\MM$ is called \emph{convergent} when there is a unique terminal set $J$. In that case, $\vv_J$ is the
unique stationary distribution for $\MM$ and with probability $1$, the time average
for $f : I \to \R$ converges to the $\vv_J$ expected value regardless of the initial position.

Suppose $X$ and $Y$ play Markov plans leading to the $4 \times 4$ matrix $\MM$ of (\ref{5})
and
$\vv_J$ is the stationary distribution for a terminal set $ J \subset \{ cc, cd, dc, dd \}$. If the sequence of
outcomes enters $J$, then it remains in $J$ and
with probability $1$ the average payoff   $ (s^N_X,s^N_Y) $ converges with
\begin{equation}\label{5a}
\lim_{N \to \infty}  (s^N_X,s^N_Y)  \ = \ v_1 (R,R) + v_2 (S,T) + v_3 (T,S) + v_4 (P,P),
\end{equation}where $\vv_J = (v_1, v_2, v_3, v_4)$.

If there is more than one terminal set $J$, then with probability
$1$ the sequence will enter some terminal set $J$ with the probabilities
for different $J$'s depending upon the initial plays.

For example, suppose that $\pp$ satisfies $0 < p_i < 1$ for $i = 1,\dots,4$ and that $\qq$ satisfies the analogous condition.
Every entry of the associated Markov matrix $\MM$ is positive and $\{ cc, cd, dc, dd \}$ is the unique terminal set.
So there is a unique stationary distribution $\vv$ with $v_i > 0$ for $i = 1,\dots,4$ and with probability one the
outcome sequence passes repeatedly through each of the four outcomes with the average payoff sequence
converging according to (\ref{5a}).
\vspace{.5cm}

Smale in \cite{Sm} aggregates the data in a different way. He suggest using as data the current average payoff given by (\ref{4}).
A \emph{Smale plan} is a function $\pi : \S \to [0,1]$. If X uses the Smale plan $\pi$ then in round $N + 1$ he plays c with
probability $\pi(s^N_X,s^N_Y)$. Again we have the switch due to reverse in labeling for the other player. Let
$Switch : \S \to \S$ be defined by $Switch(s_X,s_Y) = (s_Y,s_X)$.  If Y uses the Smale plan $\pi$ then she cooperates
with probability $\pi \circ Switch(s^N_X,s^N_Y) = \pi(s^N_Y,s^N_X)$.  That is, Y responds to $s \in \S$ by
using $\pi \circ Switch$ applied to $s$.

So if X uses $\pi_X$ and Y uses $\pi_Y$ then the outcomes $cc, cd, dc, dd$ at time $N + 1$ have probabilities given by
\begin{align}\label{6}
\begin{split}
p_{cc} \ = \ &\pi_X(s^N) \cdot \pi_Y(Switch(s^N)), \\
p_{cd} \ = \ &\pi_X(s^N) \cdot (1 - \pi_Y(Switch(s^N))),\\
p_{dc} \ = \ &(1 - \pi_X(s^N)) \cdot \pi_Y(Switch(s^N)),\\
p_{dd} \ = \ &(1 - \pi_X(s^N)) \cdot (1 - \pi_Y(Switch(s^N))).
\end{split}
 \end{align}
After the randomization is applied, we obtain the time $N+1$ payoff $S^{N+1} = (S^{N+1}_X,S^{N+1}_Y)$.

 Smale only uses pure strategy responses
for which $\pi$ maps to $\{ 0, 1 \}$. For the most part we will follow him. We will call $\pi^{-1}(1) \subset \S$ the
\emph{cooperation zone} for $\pi$ and $\pi^{-1}(0) \subset \S$ the
\emph{defection zone}.

The following \emph{Separation Theorems} will be used repeatedly.

\begin{lem}\label{lem01} Let $L : \R^2 \to \R$ be a nonconstant affine map, i.e. $L(x,y) = ax + by +c $ with
 $(a,b) \not= (0,0)$. Let $M$ be the maximum value of $|L|$ on $\S$. X and Y use general strategies.

If there exists a time $N^*$ so that for all $N \geq N^*$, $L(s^N) > 0$ implies
 $L(S^{N+1}) \leq 0$, then for $N \geq N^*$
  \begin{equation}\label{7}
 L(s^N) \quad \leq \quad \frac{M N^*}{N}.
 \end{equation}
 So  $\lim \sup_{N \to \infty} \ L(s^N) \leq 0$.
  \end{lem}

{\bf Proof:}  Since an affine map commutes with convex combinations, (\ref{4a}) implies
  \begin{equation}\label{8}
L(s^{N+1})    =   \frac{N}{N+1}L(s^N)  +   \frac{1}{N+1} L(S^{N+1}).
\end{equation}
So if $N \geq N^*$
\begin{equation}\label{9}
L(s^N) > 0 \  \Longrightarrow \  L(s^{N+1}) \leq   \frac{N}{N+1} L(s^N).
\end{equation}
On the other hand, since $N^* \geq 1$
\begin{equation}\label{10}
L(s^N)  \leq 0 \  \Longrightarrow \  L(s^{N+1})   \leq  \frac{1}{N+1}  L(S^{N+1})  \leq  \frac{ M N^*}{N+1}.
\end{equation}
Finally, observe that $L(s^{N^*}) \leq \frac{M N^*}{N^*}$. So inequality (\ref{7}) follows from
(\ref{9}) and (\ref{10}) by mathematical induction.

$\Box$ \vspace{.5cm}

\begin{lem}\label{lem01x}  Let $L : \R^2 \to \R$ be a nonconstant affine map. X and Y use general strategies.

(a) If $L(P,P), L(T,S) < 0$, then there is a positive integer $k$ such that for all $N$ there exists
$n$ with $N \leq n \leq kN$ such that either $L(s^n) < 0$ or X plays $c$ on round $n$.

(b) If $L(R,R), L(S,T) > 0$, then there is a positive integer $k$ such that for all $N$ there exists
$n$ with $N \leq n \leq kN$ such that either $L(s^n) > 0$ or X plays $d$ on round $n$.
\end{lem}

{\bfseries Proof:} (a) Let $m = \min \{ -L(P,P), -L(T,S) \}$ and let $k$ be an  integer  with $k \geq 2$ such that
 $\frac{M}{k - 1} < m$. Assume X plays $d$ in rounds $N, \dots, kN$
 then in each round the outcome is
 either $dc$ or $dd$ and so $L(S^{n}) \leq -m$ for $n = N+1, \dots, kN$.
 On the other hand, $L(S^n) \leq M$ for $n = 1, \dots, N$. Hence
 \begin{equation}\label{8c}
 L(s^{kN})  =  \frac{1}{kN} \Sigma_{t =1}^{kN} L(S^t)  \leq  \frac{1}{kN}(NM -(k - 1)Nm) < 0.
 \end{equation}

The proof of (b) is completely analogous.

$\Box$ \vspace{.5cm}

\begin{theo}\label{theo01a} Let $L : \R^2 \to \R$ be a nonconstant
 affine map with $M$ be the maximum value of $|L|$ on $\S$. Player X uses a Smale plan $\pi$ from round $N^*$ on and
player Y uses an arbitrary plan and X and Y use arbitrary initial plays.

 (a)  Assume that  $L(s) > 0$ implies $\pi(s) = 0$, i.e. $\{ L > 0 \}$ is contained in the defection zone of $\pi$.
\begin{itemize}

\item[(i)] If  $L(P,P), L(T,S) \leq 0$ then for all $N \geq N^* $.
 \begin{equation}\label{7a}
 L(s^N) \quad \leq \quad \frac{MN^*}{N}.
 \end{equation}
 So  $\lim \sup_{N \to \infty} \ L(s^N) \leq 0$.

 \item[(ii)] If  $ L(P,P), L(T,S) < 0$ then there is a positive integer $k$ so that for all $N \geq N^*$ there exists
 $n$ with $N \leq n \leq kN$  such that $L(s^n) \leq 0$.

 \item[(iii)] If $ L(R,R), L(S,T) > 0$ then  there is a positive integer $k$ so that for all $N \geq N^*$ there exists
 $n$ with $N \leq n \leq kN$  such that
 X plays $d$ on round $n$.
 \end{itemize}

(b) Assume that $L(s) < 0$ implies $\pi(s) = 1$ i.e. $\{ L < 0 \}$ is contained in the cooperation zone of $\pi$.
\begin{itemize}

\item[(i)] If  $L(R,R), L(S,T) \geq 0$ then for all $N \geq N^* $.
 \begin{equation}\label{7b}
 L(s^N) \quad \geq \quad -\frac{M N^*}{N}.
 \end{equation}
 So  $\lim \inf_{N \to \infty} \ L(s^N) \geq 0$.

 \item[(ii)] If  $ L(R,R), L(S,T) > 0$ then there is a positive integer $k$ so that for all $N \geq N^*$ there exists
 $n$ with $N \leq n \leq kN$  such that $L(s^n) \geq 0$.

 \item[(iii)] If $ L(P,P), L(T,S) > 0$ then  there is a positive integer $k$ so that for all $N \geq N^*$ there exists
 $n$ with $N \leq n \leq kN$  such that
 X plays $c$ on round $n$.
   \end{itemize}
\end{theo}

 {\bf Proof:} (a)(i)  If $L(s^N) > 0$ then $X$ plays $d$ and so
 the $N + 1$ outcome is either $dc$ or $dd$. Hence, $S^{N+1}$ is either $(T,S)$ or $(P,P)$ which implies $L(S^{N+1}) \leq 0$.
So we can apply Lemma \ref{lem01}    to get (\ref{7a}).

(a)(ii) Apply Lemma \ref{lem01x}(b) to obtain $k$ so that for some $n$ between $N$ and $kN$ either $L(s^n) < 0$ or
X plays $c$ on round $n$. By assumption, if X plays $c$ on round $n$ then $L(s^n) \leq 0$.

(a)(iii) Apply Lemma \ref{lem01x}(a) to obtain $k$ so that for some $n$ between $N$ and $kN$ either $L(s^n) > 0$ or
X plays $d$ on round $n$. If $L(s^n) > 0$ then X plays $d$ on round $n$.

 The proof of (b) is completely analogous to that of (a) applying Lemma \ref{lem01} to $-L$
to get (\ref{7b}).

$\Box$ \vspace{.5cm}

We will say that a player \emph{eventually} uses a plan when there exists $N^*$ so that the plan is used for all $N \geq N^*$.

{\bf Notation:} For distinct points $A, B \in \R^2$ we will use $[A,B]$ for the closed segment connecting the points,
$[A,B)$ for the half-open segment $[A,B] \setminus B$, etc.  We will denote by $)A,B($ the line through
$A$ and $B$ and use $[A,B($ for the ray from $A$ through $B$. In general, for a finite set of points $\{A_1,\dots,A_n \}$
we will use $[A_1,\dots,A_n ]$ for the convex hull. We will call the line $) (P,P), (R,R) ($ the \emph{diagonal} and
the line $) (S,T), (T,S) ($ the \emph{co-diagonal}. If $\ell_1$ and $\ell_2$ are non-parallel lines, then we
will let $\ell_1 \cap \ell_2$ denote the point of intersection, abusively identifying the singleton set with the point
contained therein.
\vspace{.5cm}

Any line $\ell$ in $\R^2$ is the intersection of two half-planes $H^+$ and $H^-$. When the line is not vertical
we will use $H^+$ for the upper half-plane. We will then refer to $H^+ \setminus \ell$ as the set of points above $\ell$,
with $H^- \setminus \ell$ the points below $\ell$. Up to multiplication by a positive constant
an affine map $L : \R^2 \to \R$ is uniquely defined by the conditions that $L$ is zero on $\ell$ and is positive on
$H^+ \setminus \ell$. We will say that $L$ is \emph{associated with} $\ell$ and vice-versa.

A line $\ell$ is called a \emph{separation line} for the game when $(S,T)$ and $(R,R)$ lie in one half-plane while
$(P,P)$ and $(T,S)$ lie in the other.

A separation line intersects the segments $[(S,T), (P,P)]$ and\\ $[(R,R),(T,S)]$ and so is determined by a
choice of a point on each of these segments.  Any point on $[(S,T), (P,P)]$ may be used.  If $P \leq \frac{1}{2}(T + S)$
then any point on $[(R,R),(T,S)]$ may be used. However,  if $P > \frac{1}{2}(T + S)$ then the line
 $)(S,T), (P,P)($ intersects the open segment $((R,R),(T,S))$ at a point we will label $\bar W$.
The point $\bar W$ is on or below any separation line.
In general, a separation line has slope $m$ with $|m| \leq 1$, with $m = 1$ only
for the diagonal line   and with $m = -1 $ only for the co-diagonal.
The co-diagonal is a separation line only when $P \leq \frac{1}{2}(T + S)$.

A separation line cannot be vertical, so $(S,T), (R,R) \in H^+$ and
$(P,P), (T,S) \in H^-$. Hence, for an associated affine map $L$,
\begin{equation}\label{11a}
L(R,R), L(S,T) \ \geq \ 0 \ \geq \ L(P,P), L(T,S). \hspace{2cm}
\end{equation}

\begin{cor}\label{cor01e} Assume that eventually X plays  a Smale plan $\pi$, Y uses an arbitrary plan and that
the initial play is arbitrary.  Let $\Omega$ be the limit point set of an associated sequence of outcomes. Let $C \subset \S$
be a closed, convex set and $\ell$ be a separation line.

(a) If $(P,P), (T,S) \in C$ and $\S \setminus C$ is contained in the defection zone $\pi^{-1}(0)$, then $\Omega \subset C$.

(b) If $(R,R),(S,T)  \in C$ and $\S \setminus C$ is contained in the cooperation zone $\pi^{-1}(1)$, then $\Omega\subset C$.

(c) If $\S \cap \ell \subset C$ and $\pi(s) = 0$ for $s$ above $C$ and $\pi(s) = 1$ for $s$ below $C$, then $\Omega \subset C$.
\end{cor}

{\bf Proof:} (a): For $s \in \S \setminus C$, let $s'$ be the closest point in $C$. The line $\ell$ through $s'$ which is
perpendicular to $)s,s'($ is a line of support for $C$.  That is, if $L$ is an affine function associated with $\ell$
such that $L(s) > 0$ then $L \leq 0$ on $C$. From Theorem \ref{theo01a} (a)(i) it follows that $L \leq 0$ on $\Omega$.
In particular, $s \not\in \Omega$.

(b): Proceed as above, using Theorem \ref{theo01a} (b)(i), instead.

(c): Let $C^+$ consist of the points of $\S$ on  or above $C$ and $C^-$ consist of the points of $\S$ on or below $C$.
To be precise, if $H^{\pm}$ are the half-spaces associated with $\ell$, $C^{\pm} = C \cup (H^{\pm} \cap \S)$.
These are each closed, convex sets. Because $\ell$ is a separation line, $(P,P), (T,S) \in C^-$ and
$(R,R),(S,T)  \in C^+$. From (a) it follows that $\Omega \subset C^-$ and from (b) that $\Omega \subset C^+$. Thus,
$\Omega \subset C^- \cap C^+ = C$.

$\Box$ \vspace{.5cm}

\begin{df}\label{def02} The map $\pi : \S \to [0,1]$ is a \emph{simple Smale plan} with separation line $\ell$ if
for an associated affine function $L$ for $\ell$
\begin{equation}\label{11}
\pi(s) \quad = \quad \begin{cases} 0 \quad \text{if} \ L(s) > 0, \\
1 \quad \text{if} \ L(s) < 0.\end{cases}
\end{equation}
\end{df}
\vspace{.5cm}

That is, the $\pi$ player responds with $d$ if $s$ is above the line $\ell$ and with $c$ if $s$ is below the line.
We do not specify the value of $\pi$ on the line.

\begin{cor}\label{cor01b} Assume that from some round $N^*$ on, player X uses a simple Smale plan with separation line
$\ell$ and associated affine function $L$. Let $M$ be the maximum  of $|L|$ on
$\S$. Assume that
player Y uses an arbitrary plan and that the initial plays are arbitrary.

For all $N \geq N^*$
 \begin{equation}\label{7c}
 |L(s^N)| \quad \leq \quad \frac{MN^*}{N}.
 \end{equation}
 So $\lim _{N \to \infty} \ L(s^N) = 0$ and the limit point set $\Omega$ is contained in $\ell$.
 \end{cor}

 {\bf Proof:} Clearly, (\ref{7c})
 follows  from (\ref{7a}) and (\ref{7b}), given (\ref{11a}) and (\ref{11}).

$\Box$ \vspace{.5cm}

Thus, if eventually X plays a simple Smale plan and player Y uses an arbitrary plan,
then after the randomization for mixed strategies
has been applied, a sequence of outcomes is obtained and
the limit point set $\Omega$ of the corresponding payoff sequence $\{ s^N \}$ is a point or closed segment in the
separation line $\ell$ by Proposition \ref{prop00}.

The plan All-C, which always cooperates, and so has $\pi = 1$ on $\S$, is a simple
Smale plan with separation line  $)(R,R),(S,T)($.
If $P \leq \frac{1}{2}(T + S)$ then All-D, which always defects,  and so has $\pi = 0$ on $\S$, is a simple
Smale plan with with separation line  $)(P,P),(T,S)($. However, if $P > \frac{1}{2}(T + S)$
then every simple Smale plan cooperates below the line $)(P,P),(T,S)($ and so, in particular,
has $\pi = 1$ on a neighborhood of $\frac{1}{2}(T + S,T + S)$.

If $P \leq E \leq R$ then the horizontal line $\{ s_Y = E \}$ is a separation line. The associated simple Smale
plan is the Smale version of an \emph{equalizer plan} introduced in \cite{BNS}   and also described by Press and Dyson \cite{PD}.
If X uses this equalizer plan then the limiting payoff for Y is E regardless of Y's play. On the other hand, the payoff
to X can be anything between R and P, or even lower.

\vspace{.5cm}
\includegraphics{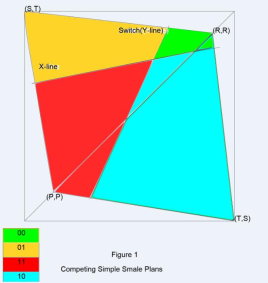}
\vspace{.5cm}

If, eventually, X and Y play simple Smale plans $\pi_X$ and $\pi_Y$ with separation lines $\ell_X$ and $\ell_Y$, respectively,
then Y responds with $\pi_Y \circ Switch$ and so the set $\Omega$ of limiting payoffs lies on the intersection
$\ell_X \cap Switch(\ell_Y)$. Except for the extreme cases with $\ell_X = \ell_Y = Switch(\ell_Y) $ equal to the diagonal
or $\ell_X = \ell_Y = Switch(\ell_Y) $ equal to the co-diagonal the intersection is a single point and so
the payoff sequence $\{ s^N \}$ converges to the intersection point $\ell_X \cap Switch(\ell_Y)$.

\begin{prop}\label{prop01c} (a) If, eventually, X plays the plan All-C then for any strategy for Y and any initial plays,
the limit point set $\Omega$ is contained in the segment $[(R,R),(S,T)]$.

(b) If, eventually, X plays the plan All-D then for any strategy for Y and any initial plays,
the limit point set $\Omega$ is contained in the segment $[(P,P),(T,S)]$.
\end{prop}

{\bfseries Proof:} (a) Since All-C is a simple Smale plan with separation line $\ell_1  \ = \ )(R,R),(S,T)($ the limit point set
lies in $\ell_1 \cap \S = [(R,R),(S,T)]$. The result also follows from Corollary \ref{cor01e} (b) with $C = [(R,R),(S,T)]$

(b) If $P \leq \frac{1}{2}(T + S)$ then All-D is a simple Smale plan with separation line $\ell_2  \ = \ )(P,P),(T,S)($ and
we can proceed as in (a). If $P > \frac{1}{2}(T + S)$, All-D is not a simple Smale plan.  The result nonetheless follows
 from Corollary \ref{cor01e} (a) with $C = [(P,P),(T,S)]$.

$\Box$ \vspace{.5cm}

\begin{cor}\label{cor01d} Assume that, eventually, X uses a Smale plan $\pi$. Let Y use an arbitrary strategy and let X and Y
use arbitrary initial plays.  If $\Omega$ is the limit point set, then
$\Omega$ is not contained in the $\S$ interior of $\pi^{-1}(1) \setminus [(R,R),(S,T)]$ and
$\Omega$ is not contained in the $\S$ interior of $\pi^{-1}(0) \setminus [(P,P),(T,S)]$.
\end{cor}

{\bf Proof:} Assume X uses $\pi$ from some time $N^*_1$.
Let $U$ denote the interior of $\pi^{-1}(1) \setminus [(R,R),(S,T)]$. If $\Omega \subset U$ then by
Proposition \ref{prop00} there exists $N^* \geq N_1^*$ such that $s^N \in U$ for all $N \geq N^*$. This implies that for
every round beyond $N^*$, X plays $c$. So the sequence of payoffs is the same as though X plays All-C
starting from round $N^*$. So by Proposition \ref{prop01c} (a) the limit point set would be contained in  $[(R,R),(S,T)]$.
This contradiction shows that $\Omega \subset U$ is impossible.

The second assertion similarly follows from Proposition \ref{prop01c} (b).

$\Box$
\vspace{1cm}

\section{Good Plans}
\vspace{.5cm}

We describe informally the conditions that we would like a \emph{good plan} to satisfy.

\begin{itemize}

\item \textbf{(Cooperation Condition)} If the players X and Y use  fixed good strategies,
i.e. good plans together with initial cooperation, then for every $N$,
$S^N \ = \ (R,R)$ and so, of course, the time averages $s^N \ = \ (R,R)$ for all $N$ as well.

\item \textbf{(Protection Condition)} If X eventually plays a fixed good plan,
and $s^* = (s_X^*,s_Y^*)$ is a limit point for the sequence
$\{ s^N \}$ with arbitrary initial play and with Y using any plan, then
\begin{equation}\label{12}
s_Y^* \geq R \qquad \Longrightarrow \qquad s_X^*  \ = \ s_Y^*  \ = \  R.
\end{equation}

\item \textbf{(Robustness Condition)} If eventually the players X and Y use
fixed good plans, then regardless of the earlier play \\
$\lim_{N \to \infty} \ s^N \ = \ (R,R)$ at least with probability one.

\end{itemize}

A Markov plan is called \emph{agreeable} if the response to a $cc$ outcome is always $c$. That is,
$\pp $ satisfies $ p_1 = p_{cc} = 1$. A Markov plan
is called \emph{firm} if the response to a $dd$ outcome is always $d$. That is,
$\pp $ satisfies $ p_4 =  p_{dd} = 0$. For example, the $TFT$ plan with $\pp = (1, 0, 1, 0)$ is both agreeable and firm.

If both X and Y use Markov plans then $\{ cc \}$ is a terminal set if and only if both plans are agreeable.

We call a general plan is \emph{weakly agreeable} if the response is $c$ when every previous outcome is $cc$.
A  general plan is called \emph{weakly firm} if the response is $d$ when every previous outcome is $dd$.
For example, a Smale plan $\pi$ is weakly agreeable if and only if $\pi(R,R) = 1$ and is weakly firm if
and only if $\pi(P,P) = 0$.
Clearly,  a Markov plan is weakly agreeable if and only if it is agreeable and is weakly firm if and only if it is firm.

If X and Y both use weakly agreeable plans and initially cooperate then every outcome is $ cc $ and
$s^N = (R,R)$ for all $N$.

So to obtain the Cooperation Condition we demand that a good plan be weakly agreeable.

The agreeable Markov plans which satisfy the Protection Condition can be completely characterized.
\vspace{.5cm}

\begin{theo}\label{theo03a} Let $\pp = (p_1, p_2, p_3, p_4)$ be an agreeable Markov plan so that $p_1 = 1$.

The plan $\pp$ satisfies the Protection Condition if and only if   the following inequalities hold.
\begin{equation}\label{ineqa}
\frac{T - R}{R - S} \cdot p_3 \ < \ (1 - p_2) \qquad \mbox{and} \qquad  \frac{T - R}{R - P} \cdot p_4 \ < \ (1 - p_2).
\end{equation}
\end{theo}
\vspace{.5cm}

{\bf Proof:}  See  \cite{A-16} Theorem 1.5, where a plan is called good if it is agreeable and satisfies the Protection Condition.
See also \cite{A-13}.

$\Box$ \vspace{.5cm}

{\bf Remark:} Note that (\ref{ineqa}) implies $p_2 < 1$.
\vspace{.5cm}

For a Markov plan $\pp$ we will refer to (\ref{ineqa}), together with the assumption $p_1 = 1$,
as the \emph{protection inequalities}.
Thus, a Markov plan is agreeable and  satisfies the  Protection Condition exactly when the the protection inequalities hold.

For Smale plans we have

\begin{theo} \label{theo03b}  Let $\ell$
be a separation line with associated affine function $L$ such that
\begin{equation}\label{ineqb}
L(R,R) = 0, \quad \text{and}  \quad L(P,R) > 0 .
\end{equation}
That is, $\ell $ is a line through  $(R,R)$ with slope $m$ satisfying $0 < m \leq 1$.

Let $\pi : \S \to [0,1]$ be a Smale plan.
If $L(s) > 0 $ implies $\pi(s) = 0$, i.e. $\{ L > 0 \}$ is contained in the defection zone,
then $\pi$ satisfies the Protection Condition.

In particular, if $\pi$ is a simple Smale plan with separation line $\ell$ then $\pi$ satisfies the Protection Condition.
\end{theo}

{\bf Proof:} The line $\ell$ contains $(R,R)$, and  the point $(P,R)$ lies above $\ell$. Since $\ell$ is a
separation line, it follows that $s = (R,R)$ is the only point of $\S \cap H^-$ with $s_Y \geq R$.
 By Theorem \ref{theo01a} (a) $L(s^*_X,s^*_Y) \leq 0$ for every limit point $s^*$, i.e. $\Omega \subset \S \cap H^-$.

Hence, $(R,R)$ is the only possible limit point $s^*$ with $s^*_Y \geq R$.

$\Box$ \vspace{.5cm}

If $\ell$  is a line through $(R,R)$ with slope $m$ satisfying $0 < m \leq 1$ then $\ell$ is a separation line.
and we call such a line $\ell$ a \emph{protection line} for $\pi$ if $\pi = 0$ above the line.
 The above result says exactly that if a Smale plan admits a protection line then it satisfies the Protection Condition.

\begin{df}\label{def04a} A Markov plan $\pp = (p_1, p_2, p_3, p_4)$ is \emph{generous}
when
\begin{equation}\label{13}
p_1 \ = \ 1, \quad 1 > p_2 > 0, \quad p_4 > 0.
\end{equation}
\end{df}
\vspace{.5cm}

That is, a generous plan is agreeable and with positive probability responds to an opponent's defection with cooperation,
but does not always cooperate from a $cd$ outcome (which had payoff $(S,T)$).

\begin{theo}\label{theo04a} Assume that X and Y, eventually, play Markov plans $\pp$ and $\qq$ respectively.
If both plans are generous then $\{ cc \}$ is the only terminal set for the associated Markov chain.
So from any initial play, with probability one, eventually the
outcome sequence is constant at $cc$ and so $\lim_{N \to \infty} \ s^N \ = \ (R,R)$. \end{theo}

{\bf Proof:} Since the two plans are agreeable, $\{ cc \}$ is a terminal set.  Since $p_4, q_4 > 0$, there is a
positive probability of moving from $dd$ to $cc$ and so $dd$ is a transient state. From $cd$, $1 > p_2 > 0$ implies
that X plays either $c$ or $d$ with positive probability. If Y plays $c$ (or $d$) with positive probability
then from $cd$ there is a positive probability of moving to $cc$ (resp.to $dd$ and thence to $cc$).
Hence, $cd$ is transient.  Recall that Y uses the response vector $(q_1, q_3, q_2, q_4)$ and so cooperates after
$dc$ with probability $q_2$. Thus, a symmetric argument shows that $dc$ is transient.

If the players use $\pp$ and $\qq$ from time $N^*$ on then from that point, the play follows the Markov chain given by
$\MM$ and so with probability one the sequence of outcomes eventually arrives at the unique terminal set $\{ cc \}$.

$\Box$ \vspace{.5cm}

\begin{df}\label{def04b} A Smale plan $\pi : \S \to [0,1]$ is \emph{generous}
if $\pi(s) = 1$ when $s_X \geq s_Y$ and there exists an open subset $U$ of $\S$
which contains the half-open segment $ [\frac{1}{2}(T+S,T+S), (R,R)) $
such that  $\pi(s) = 1$ for $s \in U$. That is, the cooperation zone contains $U$ and
the points on and below the diagonal line.
\end{df}
\vspace{.5cm}

The following is essentially a part of \cite{Sm} Theorem 1.

\begin{theo}\label{theo04b} Assume  that, eventually, X and Y  play  Smale plans $\pi_X$ and $\pi_Y$, respectively.
If both plans are generous then  from any initial state in $\S$,
$\lim_{N \to \infty} \ (s_X^N,s_Y^N) \ = \ (R,R)$. \end{theo}

{\bf Proof:} Assume that both plans are adopted by time $N^*$.
Let $L_0(s) = s_Y - s_X$. That is, $L_0$ is the affine function associated with the diagonal separation line.
Let $L_1(s) = s_X + s_Y - T - S$. That is, $L_1$ is the affine function associated with the co-diagonal.
Notice that the maximum value of $L_1$ on $\S$ is $L_1(R,R) = 2R - T - S$ and that $L_1 < L_1(R,R)$ on $\S \setminus \{ (R,R) \}$.

Since $L_0(s) < 0$ implies $\pi_X(s) = 1$ and $\pi_Y(s) = 1$, we can apply Theorem \ref{theo01a} (b) to $L_0$ and $\pi_X$ to get
$L_0(s^N) \geq - M_0N^*/N$,  where $M_0 = T - S$ is the maximum value of $|L_0|$ on $\S$.
The Y player uses $\pi_Y \circ Switch$ and so we apply the theorem to $L_0 \circ Switch$  and $\pi_Y \circ Switch$ to get
$-L_0(s^N) = L_0 \circ Switch(s^N) \geq -M_0N^*/N$. Hence,
\begin{equation}\label{14}
|L_0(s^N)| \quad \leq \quad M_0N^*/N. \hspace{1cm}
\end{equation}
Thus, the limit point set $\Omega$ is contained in the diagonal.

Observe  that everywhere in $\S$ either  X or Y plays $c$ and so the only outcomes
are $cd, dc$ and $cc$. Hence,  we  have $L_1(S^N) \geq 0$ for all $N$.
We can apply Lemma \ref{lem01} with  $L = - L_1$ to get $L_1(s^N) \geq - \frac{M_1 N^*}{N}$ for all $N \geq N^*$,
where $M_1$ is the maximum value of $|L_1|$ on $\S$.
Hence, on $\Omega$, $L_1 \geq 0$.
Notice that in the case when $P \geq (T + S)/2$, $L_1 \geq 0$ on all of $\S$.

Thus, $\Omega \subset [\frac{1}{2}(T + S, T + S), (R, R)]$.

Let $U = U_X \cap Switch(U_Y)$ where $U_X$ and $U_Y$ are the open sets containing the half-open interval
$[\frac{1}{2}(T + S,T + S),(R,R))$ on which   $\pi_X = 1$ and $\pi_Y = 1$, respectively. By assumption, if
$s^N \in U$ then both players play $c$ with outcome $cc$ at time $N+1$. So $S^{N+1} = (R,R)$.

For
$ \e > 0,  $ $\{ L_1 > L_1(R,R) - \e \}$ is a neighborhood of $(R,R)$ and so $U \cup \{ L_1 > L_1(R,R) - \e \}$
is a neighborhood of   $[\frac{1}{2}(T + S, T + S), (R, R)]$ and so of $\Omega$. From  Proposition \ref{prop00}
it follows that there exists $N_{\e} \geq N^*$ so that $N \geq N_{\e}$ implies
$s^N \in U \cup \{ L_1 > L_1(R,R) - \e \}$. Hence, for $N \geq N_{\e}$, $L_1(s^N) - L_1(R,R) + \e \leq 0$
implies $s^N \in U$ and so $L_1(S^{N+1}) = L_1(R,R)$ and so $L_1(S^{N+1}) - L_1(R,R) + \e = \e > 0 $.
We can again apply Lemma \ref{lem01} this time to $L = - (L_1 - L_1(R,R) + \e)$ to get
\begin{equation}\label{15a}
N \geq N_{\e} \quad \Longrightarrow \quad L_1(s^N) \quad \geq \quad L_1(R,R) - \e - \frac{(M + \e)N_{\e}}{N},
\end{equation}
where $M$ is the maximum of $|L_1 - L_1(R,R)|$ on $\S$.  Hence, there exists $N'_{\e} \geq N_{\e}$ so that
\begin{equation}\label{16a}
N \geq N'_{\e} \quad \Longrightarrow \quad L_1(s^N) \quad \geq \quad L_1(R,R) - 2\e.
\end{equation}

Thus, $\lim_{N \to \infty} \ L_1(s^N) = L_1(R,R)$.

Since, $(R,R)$ is the unique maximum point for $L_1$,
it follows that $\lim_{N \to \infty} \  s^N  \ = \ (R,R)$.

$\Box$ \vspace{.5cm}

It would be nice to show that if X plays a generous Smale plan and Y plays a generous Markov
plan, then with probability one $\lim_{N \to \infty} \  s^N  \ = \ (R,R)$. I don't know if the conjecture
is true in full generality. However, we get the result we want if we
strengthen the assumption.

\begin{df}\label{def04c} A Smale plan $\pi$ is \emph{convex-generous} if $\pi : \S \to \{ 0,1 \}$,i.e
no mixed strategy responses, and the cooperation zone $\pi^{-1}(1)$ is a closed convex set $C$ such that

\begin{itemize}
\item[(i)] $(P,P), (R,R), (T,S) \in C$.

\item[(ii)] $(S,T) \not\in C$.

\item[(iii)]  $\frac{1}{2}(T + S,T + S) \in C^{\circ}$
where $C^{\circ}$ is the interior of $C$ with respect to $\S$.
\end{itemize}
\end{df}
\vspace{.5cm}

By (i) and (iii), there exists $t^* \in (\frac{1}{2},1)$ such that
 $(1-t)(T,S) + t(S,T) \in C$ if and only if $0 \leq t \leq t^*$. Let $V$ denote the point
$(1-t^*)(T,S) + t^*(S,T)$.

Let $\bar P = \min(P, \frac{1}{2}(T + S))$. Thus, the diagonal line  intersects $\S$ in the segment
$[(R,R),(\bar P, \bar P)]$ and $(\bar P, \bar P) \in C$, by (i) and (iii). Hence,
 the convex hull of $[ V, (\bar P, \bar P), (R,R),(T,S) ]$ is
contained in $C$. This contains $(s_X,s_Y) $ if $s_X \geq s_Y$. Furthermore, its $\S$ interior contains
$[(\frac{1}{2}(T + S,T + S), (R,R))$ and so, as expected, $\pi$ is generous.
\vspace{.5cm}

\begin{theo}\label{theo05} Assume X eventually  plays a convex-generous
Smale plan and Y plays a generous Markov plan. With
probability one there is a time after which both players play $c$ and so
$\lim_{N \to \infty} \ s^N \ = \ (R,R)$. \end{theo}

{\bf Proof:} Let X play $\pi$ with convex set $C = \pi^{-1}(1)$ and let Y use $\qq = (q_1,q_2,q_3,q_4)$ with
$q_1 = 1$ and $\e < q_2,  \e < q_4$ for some $\e  > 0$. Recall that for Y, $q_2$ is the probability of cooperating in a round
following the outcome $dc$.

{\bf Claim 1:} With probability one, for every $N$ there exists $n \geq N$ such that $s^n \in C$.

Proof: It suffices to show that for every $N$ the event
$$E_N = \{ s^n \in \S \setminus C : \ \text{ for all} \ n \geq N \}$$
has probability zero.

Let $\ell_1  \ = \ )V,(R,R)($, which is a separation line, and let $L_1$ be an affine function associated with
$\ell_1$. Clearly, $(P,P)$ and $(T,S)$ lie below the line and so  $L_1(P,P), L_1(T,S) < 0$.
By Lemma \ref{lem01x} (a) there exists
a positive constant $k_1$ (which depends only on $L_1$)
such that for some $N_1$ between $N$ and $k_1 N$ either $L_1(s^{N_1}) < 0$ or X plays $c$ on round $N_1$.
Note that the latter is equivalent to $s^{N_1} \in C$.

Case 1 [$P \geq \frac{1}{2}(T+S)$]: The point $V$ lies on the side $[(S,T),(T,S)]$ of the triangle $\S$. Thus,
every point of $\S$ on or below $\ell$ lies in $C$. Hence,
$L_1(s) \leq 0$ implies $s \in C$ and so $s^{N_1} \in C$ in any case.

Case 2 [$P < \frac{1}{2}(T+S)$]: The line $\ell_1$ intersects the open segment $((P,P),(S,T))$ in a point $V'$.
Let $\Delta$ denote the triangle $[(P,P),V',V]$ and
$\Delta' = \Delta \setminus C $. Observe that $C \cup \Delta'$ contains
the set $\{ L_1 \leq 0 \} \cap \S$. Assuming $E_N$ , $s^{N_1}$ lies in $\Delta' $ since
it is not in $C$. For $n \geq N$ if $ s^n \in \Delta'$ then,
since it is not in $C$, the payoff $S^{n+1}$ is either
$(P,P)$ or $(T,S)$  and so by (\ref{4a}) $L_1(s^{n+1}) \leq \frac{n}{n+1} L_1(s^n) \leq 0$.
Since $s^{n+1} \not\in C$, it follows that
$s^{n+1} \in \Delta' $.
Inductively, we have $s^n \in \Delta' $ for all $n \geq N_1$ and so for all $n \geq k_1 N$.

Now if among the rounds $k_1 N, \cdots, M-1$ Y plays $c$ exactly $r$ times then
\begin{equation}\label{16}
s^{k_1N + M} = \frac{ k_1 N s^{k_1N} + (M - r)(P,P) + r(T,S)}{k_1 N \ + \ M},
\end{equation}
and so the vector from $(P,P)$ to $s^{k_1N + M}$ is the nonzero vector
\begin{equation}\label{17}
s^{k_1N + M} - (P,P)= \frac{ k_1 N [s^{k_1N} - (P,P)] + r[(T,S) - (P,P)]}{k_1 N \ + \ M},
\end{equation}
Normalize these vectors to obtain $[s^{k_1N + M} - (P,P)]_1$ of length $1$.

Suppose that as $M \to \infty$ the number $r$ of $c$ plays by Y is unbounded. Then, as $M \to \infty$ these
unit vectors have $[(T,S) - (P,P)]_1$ as a limit point. On the other hand, the closed set of vectors\\
$\{ [s - (P,P)]_1 : s \in \Delta \setminus \{(P,P)\} \} = \{ [s - (P,P)]_1 : s \in [V,V'] \} $
does not contain $[(T,S) - (P,P)]_1$.

It follows that, assuming $E_N$, there exists $R < \infty$ such that Y plays $c$ at most $R$ times
and so from some $N_2$ onward Y always plays $d$.
For each $N_2 \geq k_1 N$ the event $E_{N,N_2} = E_N$ \emph{and Y plays $d$ on every round $n$ with $n \geq N_2$}
has probability zero because at each such round Y is responding to $dd$ by playing $d$. These are independent events
each with probability at most $1 - \e$. So $E_N = \bigcup_{N_2 \geq k_1 N} \ E_{N,N_2}$ has probability zero.

This completes the proof of Claim 1.

If Y plays c at any time $N$ when $s^N \in C$ then the outcome of the $N+1$ round is $cc$ with payoff $(R,R)$.
Since $C$ is convex $s^{N+1} \in C$ and so X next plays C and Y next plays C because $\qq$ is agreeable.
Inductively $cc$ is the outcome and $s^n \in C$ for every round $n$
with $n \geq N$. Let $E_0 $ denote the event \emph{Y plays $d$ whenever $s^n \in C$} and
$\tilde E = E_0 \setminus \bigcup_N \{ E_N \}$. From Claim 1, it suffices to show that $\tilde E$ has probability zero.

{\bf Claim 2:} Assuming $\tilde E$, for every $N$ there exists $n \geq N$ such that $s^n \in \S \setminus C$.

Proof: If we assume $E_0$, then whenever X plays $c$, Y plays $d$ leading to payoff $(S,T)$.
If for some $N$, $s^n \in C$ for all
$n \geq N$ and  $E_0$ is true, then for all $n > N$ we have $S^n = (S,T)$ and the sequence $\{ s^n \}$ would
converge to $(S,T)$, but the complement of $C$ is a neighborhood of $(S,T)$ and so eventually $s^n \not\in C$.
The contradiction proves Claim 2.

Assuming $\tilde E$,
$s^n \in C$ and $s^n \in \S \setminus C$ each occur infinitely often by Claim 1 and Claim 2.

Now let $N_k$ be the $k^{th}$ return time to $C$ from $X \setminus C$. This is an infinite sequence of Markov times.
Since at time $N_k - 1$ the payoff sequence was in $\S \setminus C$, X played $d$ and so at time $N_k$ Y plays $d$
in response to either a $dc$ or a $dd$. Playing $d$ in these cases has probability at most $1 - \e$. Furthermore,
the Y play at time $N_k$ is independent of the previous plays.
Thus, again $\tilde E$ requires an infinite sequence of independent
events, each with probability at most $1 - \e$. Hence, $\tilde E$ has probability zero.

$\Box$ \vspace{.5cm}

{\bf Remark:} Our assumptions on $\qq$ allow the possibility $q_3 = 0$. So if $(S,T)$ were in $C^{\circ}$ then
from a neighborhood of $(S,T)$ it would be a limit point following an infinite sequence of $cd$ outcomes.
Notice that the assumption $(S,T) \not\in C$ is analogous to the assumption that $p_2 < 1$ for a
 generous Markov plan. If $p_2 = 1, p_3 = 0$ and $\qq$ satisfies the analogous condition then
$\{cd \}$ and $\{dc \}$ are terminal sets when  X plays $\pp$ and Y plays $\qq$.
\vspace{.5cm}

\begin{df}\label{def06} We call a Markov plan  \emph{good} when it satisfies the protection inequalities (and so is
agreeable) and is generous.

We call a Smale plan  \emph{good} (or \emph{convex-good}) when it is weakly agreeable, admits a   protection line and is generous
(resp. and is convex-generous). \end{df}
\vspace{.5cm}

For example, if $\pi$ is a simple Smale plan with separation line $\ell$ then $\pi$ is good if and only if it is
weakly agreeable (i.e. $\pi(R,R) = 1$) and $\ell$ is a line through
$(R,R)$ with slope strictly between $0$ and $1$, so that $\ell$ is a protection line.
It is convex-good if and only if, in addition, $\pi = 1$ on
$\ell \cap \S$. Such a good simple Smale plan is the Smale version of what is called in \cite{HNS} a \emph{complier
strategy} (also called a \emph{generous zero-determinant strategy}, as in e.g. \cite{SP2}).
Any limit point $s^*$, other than $(R,R)$ on the separation line $\ell$ has
$s^*_Y $ larger than $s^*_X$, albeit less than $R$.

On the other hand, if $\ell$ is a line through $(P,P)$ with slope between $0$ and $1$ then it is a
separation line and the associated simple Smale plan is, when $P \leq \frac{1}{2}(T + S)$, the Smale version of what
Press and Dyson \cite{PD} call an \emph{extortionate strategy}.   Any limit point $s^*$, other than $(P,P)$ on the separation line
$\ell$ satisfies $s^*_Y < s^*_X$. Thus, if Y plays to avoid the $(P,P)$ payoff she always obtains less than X
does from the change in policy. The best reply to such an extortionate strategy is All-C. The payoff point
is then the intersection point $B = \ell \cap ((R,R),(T,S))$ with $P < B_Y < R < B_X$.

In \cite{A-16} a Markov plan, there called a memory-one plan,
is called good when it satisfies the protection inequalities, or, equivalently,
it is agreeable and satisfies the Protection Condition. The TFT plan with $\pp = (1,0,1,0)$ satisfies the
protection inequalities but is firm rather than generous. If both
X and Y use the TFT plan then from initial outcome $cc$ the sequence of outcomes is fixed at $cc$, but from
an initial $dd$, the state $dd$ is fixed. Following $cd$ or $dc$ the two states alternate leading to
convergence of the payoff sequence to $\frac{1}{2}(T+S,T+S)$.  The phenomena illustrate the failure of robustness
in the absence of generosity.

In \cite{Sm} a Smale plan $\pi$ is called good if it is generous, and so satisfies the Robustness Condition,  and
$\pi(s) = 0$ when $s_Y > R$. The line $\{ s_Y = R \}$ is an equalizer line and so Smale's conditions allows the possibility
 of a limit outcome $(s_X,R)$ with $s_X < R$.

In addition, Smale imposed the condition that $\pi(s) = 0$ when $s_X < P$. If $P > \frac{1}{2}(T + S)$, this would contradict
the condition that $\pi = 1$ when $s_X > s_Y$. Smale was only considering the case with $P < \frac{1}{2}(T + S)$ and we will
examine that situation first.

Assume $P < \frac{1}{2}(T + S)$. Choose $\ell_1$ a line through $(R,R)$ with slope strictly between $0$ and $1$, so that
the weakly agreeable simple Smale plan with separation line $\ell_1$ is good.  Let $A$ be the point of intersection
of $\ell_1$ and the open segment $((P,P),(S,T))$. Choose a point $V$ on the open segment $((R,R),A)$ with
$V_X \geq P$ or larger. Let $\ell_2$ be the line $)(P,P),V($. Let $\pi$ be a Smale plan such that $\pi(s) = 0$ if
$ s$ is above $\ell_1$ or above $\ell_2$ and $\pi(s) = 1$ at $(R,R)$, below the diagonal and on some
open set containing $[\frac{1}{2}(T+S,T+S),(R,R))$. Thus, $\pi$ is generous and since $\ell_1$ is a protection line
it is a good Smale strategy. It follows from Corollary \ref{cor01e}(c) that if X eventually plays $\pi$
against any plan for Y and any initial plays, then
the limit point set $\Omega$  is contained in the triangle  $[(P,P),(R,R),V]$.

The advantage of such a plan is that it excludes points above $\ell_2$ from the limit.  In contrast, against the good simple Smale
strategy with separation line $\ell_1$, any point of $[A,(R,R)]$ can occur as the limit if Y plays a suitable simple
Smale strategy. For example, recall that All-D is a simple Smale strategy with separation line $\ell  \ = \ )(P,P),(T,S)($.
Thus, if Y plays All-D then the limit point is the intersection point of $\ell_1$ with $Switch(\ell)$ which is $A$
with $A_Y > P > A_X$.

This sort of possible cost can always occur with a generous Markov plan $\pp$. If $Y$ plays All-D, which is a
Markov plan with $\qq = (0,0,0,0)$,
then the unique terminal set is $\{ cd,dd \}$ with stationary distribution
$\vv = (0,p_4,0,(1-p_2))/[p_4 + (1-p_2)]$. The limiting average payoff $(s_X^*,s_Y^*)$ given by (\ref{5a})
satisfies $(s_X^*,s_Y^*) - (P,P) = \e (S - P,T - P)$ with $\e = p_4/[p_4 + (1-p_2)] > 0$.  So
$s_Y^* > P > s_X^*$.

On the other hand, for plans such as $\pi$ above we have seen that the limit point set $\Omega$
against any Y play is a connected set contained in the triangle  $[(P,P),V,(R,R)]$.

However:

 \vspace{.5cm}
\includegraphics{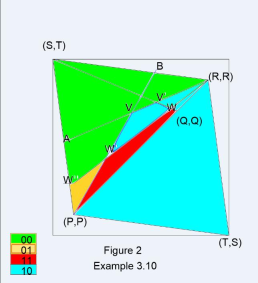}
\vspace{.1cm}

\begin{ex}  The limit set $\Omega$ need not be a point or interval and
it may contain points in the interior of the triangle.  \end{ex}

Assume $P <  \frac{1}{2}(T + S)$.  Choose $V \in \S$ with $R > V_Y > V_X \geq P$.
Let $\ell_1  \ = \ )(R,R),V($ and $\ell_2  \ = \ )(P,P),V($. Let $\pi_X$ be a Smale plan with
$\pi_X(s) = 0$ for  above $\ell_1$ or above $\ell_2$
and $\pi_X(s) = 1$ otherwise. In that case, $\pi_X$ is convex-good with $C$ the quadrilateral
$[(P,P), V, (R,R), (T,S)]$. Let $C' \subset C$ denote the triangle  $[(P,P),V,$ $(R,R)]$. As we saw above,
when X plays $\pi_X$ against any strategy for Y, the limit
point set $\Omega$  is contained in $C'$. Furthermore, by Corollary \ref{cor01d}
$\Omega$ must intersect $[(P,P),V] \cup [V,(R,R)]$.

Choose a point $V'$ on the half-open segment $[V,(R,R))$ and a point $W$ on $)(S,T),V'($ between the
lines  $\ell_1$ and the diagonal.  Let $\ell  \ = \ )(R,R),W($ so that the line $\ell$ lies between  $\ell_1$ and the diagonal,
with all three intersecting at $(R,R)$. Let $\ell' \ = \ )(S,T),W( \ = \ )(S,T),V'($.

Label the following points:
\begin{itemize}
\item $\ell_1  \ \cap  \  )(P,P),(S,T)( \ = \ A$.
\item $\ell_2  \ \cap  \   )(R,R),(S,T)( \ = \ B$.
\item $\ell'  \ \cap  \    )(P,P),(R,R)( \ = \ (Q,Q)$  \quad and \quad $ \ell' \ \cap  \ \ell_1 \ = \ V'$.
\item $\ell   \ \cap  \   \ell_2 \ = \ W'$ \quad and \quad $\ell   \ \cap  \   )(P,P),(S,T)( \ = \ W''$.
\end{itemize}

Let $\bar C$ be the quadrilateral
$[(P,P),(Q,Q),W,W'']$. Define $\pi_Y$ for Y to be the Smale plan with $\pi_Y(s) = 1$ if $s$ lies in $Switch(\bar C)$
and $= 0$ otherwise.  Recall that Y responds with $\pi_Y \circ Switch$ and so Y plays $c$ if $s \in \bar C$
and plays $d$ otherwise.

We prove that if, eventually, X plays $\pi_X$ and Y plays $\pi_Y$ then regardless of the initial play,
the limit point set $\Omega$
is the boundary of the quadrilateral $\hat L = [V,W',W,V']$ (If $V' = V, \ \Omega$ is the boundary of the
triangle  $[V,W',W]$).
\vspace{.5cm}

{\bf Proof:} Assume that X and Y play $\pi_X$ and  $\pi_Y$, respectively, beyond time $N^*$.

The lines $\ell_1, \ell_2, \ell, \ell'$ and the diagonal subdivide $\S$ into twelve polyhedral
regions. For any $\e > 0$, let
$N_{\e} > N^*$ be an $\e$ step-time after which every move from $s^N$ to $s^{N+1}$ has distance less than $\e$,
i.e. let $N_{\e}$ be greater than $\max(N^*, \sqrt{2}(T - S)/\e)$.
 Let $\e_0 > 0$ be smaller than the distance between any two non-intersecting regions and let
 $N_0 =  N_{\e_0}$.  Thus, such a small move
cannot jump between non-intersecting regions.

Let $\tilde C = \bar C \cap C$, which is the quadrilateral $[(P,P),W',W,(Q,Q)]$.

Claim: For every $N \geq N^*$ there exists $n \geq N$ such that $s^n \in \tilde C$.

First we show that the sequence of payoffs must enter $\bar C$. If not, then for every round beyond $N$,
Y plays $d$. As in the proof of Proposition \ref{cor01d}
the results after $N$ are the same as though Y plays All-D which is a simple Smale plan with
separation line $)(P,P),(T,S)($. Then $\Omega$  is contained in the intersection of the triangle $ C'$ with
$Switch( \ )(P,P),(T,S)( \ ) $  $ = \ )(P,P),(S,T)($. This intersection contains only the point $(P,P)$. If $\Omega$ were just
$(P,P)$ then for any small neighborhood $U$ of $(P,P)$ eventually $s^n \in U$. If $s^n$ is on or above the
diagonal then $s^n \in \bar C$. If $s^n$ is below the diagonal then the sequence of payoffs moves toward
$(S,T)$ and eventually makes a small jump into $\bar C$.  Either way, this contradicts the assumption that
the sequence never enters $\bar C$.

Now for $Z$ on the open segment $((P,P),W')$ let $U_Z$ consist of the points of $\S$ which are below both of
the lines $)(S,T),Z($ and $)(T,S),Z($. These are  convex open neighborhoods of $(P,P)$ which converge to
$(P,P)$ as $Z \to (P,P)$. Since $(P,P) \in \tilde C$, if $s^n \in \bar C \setminus \tilde C$ then
there exists $Z$ such that $s^n \not\in U_Z$, i.e. $s^n \in K_1 = \bar C \setminus (\tilde C \cup U_Z)$.
  Assume $n > N_0$. From such a point the sequence moves toward
$(T,S)$. During the motion it remains above $U_Z$. If the sequence does not enter $\tilde C$ from  $\bar C \setminus \tilde C$
then it jumps to below the diagonal to land in $K_2$ the set of points outside $U_Z$ which are on or below
the diagonal and $\ell'$. From such points the sequence moves back toward $(S,T)$. If it jumps over $\tilde C$ then
it re-enters $K_1$. This alternation cannot continue indefinitely.  Notice that $K_1$ and $K_2$ are a positive distance $\e_Z$
apart. Once $N \geq N_{\e_Z}$ the move from $K_1$ or $K_2$ must land in $\tilde C$.

This completes the proof of the Claim.

For any $\e > 0$, let $s^n \in \tilde C$ with $n > N_{\e}$. From this point the sequence moves toward $(R,R)$ exiting
$\bar C$ at a point above, and $\e$ close to, the line $\ell'$. The sequence now moves toward $(S,T)$, $\e$ close to and above the
line $\ell'$. It exits $C$ at a point above $\ell_1$ and $\e$ close to $V'$. Now the sequence moves toward $(P,P)$ entering
$\hat L$ on or below, and $\e$ close to $\ell_1$ and then moving back toward $(S,T)$. These $PP$ and then $ST$ alternate motions
may continue for a long time but it must eventually cease since the sequence must eventually return to $\tilde C$. The exit occurs
when the sequence lands above $\ell_2$, $\e$ close to $V$. The sequence then moves toward $(P,P)$ above and $\e$ close
to the line $\ell_2$ until it enters $\bar C \setminus \tilde C$, $\e$ close to $W'$. The sequence then moves toward
$(T,S)$ crossing $\ell_2$ close to $W'$ to re-enter $\tilde C$ now below and $\e$ close to the line $\ell$.

As $N \to \infty, \e \to 0$ and the motion gets close to motion from $W'$ to $W$, from $W$ to $V'$, from $V'$ to $V$,
and then from $V$ back to $W$.

$\Box$ \vspace{.5cm}

Turning now to the case when $P \geq \frac{1}{2}(T+S)$ we see that the additional condition imposed by Smale now does not work
so well.

Suppose you demand that $\pi$ satisfy $\pi(s) = 0$ when $s_X \leq \frac{1}{2}(T+S)$ or even just
$\pi = 0$ on the half-open segment $(\frac{1}{2}(T+S,T+S), (S,T)]$ with $\pi = 1$ on
$(\frac{1}{2}(T+S,T+S), (T,S)]$. If $\pi_X$ and $\pi_Y$ both satisfy this condition then
for $s^N \in (\frac{1}{2}(T+S,T+S), (S,T)]$ the payoff $S^{N+1} = (T,S)$
and for $s^N \in (\frac{1}{2}(T+S,T+S),(T,S)]$ the payoff $S^{N+1} = (S,T)$. Thus, the sequence remains on
$[(S,T),(T,S)]$ moving back and forth as the point $\frac{1}{2}(T+S,T+S)$ is passed with limit point
$\frac{1}{2}(T+S,T+S)$, \ul{unless} the sequence lands exactly on the point $\frac{1}{2}(T+S,T+S)$.
If that happens then, the result depends on the choices at the point $\frac{1}{2}(T+S,T+S)$.

For most initial points, this cannot happen. For example, if the initial point is an irrational mixture of
$(S,T)$ and $(T,S)$ then hitting $\frac{1}{2}(T+S,T+S)$ does not occur. However, for actual play
all outcomes are rational mixture of the four points $(S,T),(T,S),(P,P)$ and $(R,R)$. Focusing upon actual play leads to
odd results.

\begin{prop}\label{06} Suppose that $\pi(S,T) = 0, \pi(T,S) = 1$ and $\pi(s) = 1$  when $s$ lies on the diagonal.
If both X and Y play Smale plans which satisfy these conditions then for any initial plays,
$\lim s^N = (R,R)$. \end{prop}

{\bf Proof:} If the initial outcome is $cc$ or $dd$ then the initial payoff lies on the diagonal and so every successive
outcome is $cc$. If the initial outcome is $cd$, with payoff $(S,T)$ then the next outcome is
$dc$ with payoff $(T,S)$ so that $s^2 = \frac{1}{2}(T+S,T+S)$ which lies on the diagonal. Hence, in any case,
the outcome on the $n^{th}$ round is $cc$ for $n \geq 3$ and the limit result follows.

$\Box$ \vspace{.5cm}

This version of ``robustness'' is very unsatisfying. For real robustness one wants approach to $(R,R)$ even if errors occur
in the computations and if the plans are adopted only from some time $N^*$ on.

\vspace{.5cm}
\includegraphics{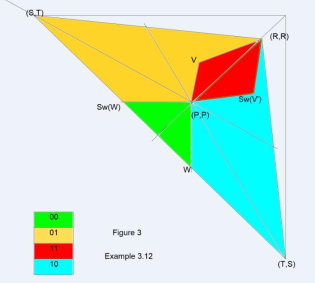}
\vspace{.1cm}

\begin{ex} Returning to the case with $P > \frac{1}{2}(T + S)$ we describe  what seems to me to be the
best version of robustness we can hope for if we demand protection against payoffs below $P$. \end{ex}

Let $\ell$ be a protection line with slope less than $1$ so that
$(P,P)$ lies below $\ell$.  Let $V \in \ell$ with $P \leq V_X < R$,
$\ell_1  \ = \ )(P,P),V($ and $\ell_2$ the vertical line $\{s_X = P \}$. Let $W$ be the point of intersection
 $\ell_2 \ \cap \ )(S,T),(T,S)($.  Assume that $\pi_X(s) = 0$ whenever $s$ is
above  $\ell $ or $\ell_1$, or if it is on or to the left of $\ell_2$. Otherwise, $\pi_X(s) = 1$.

It follows from Corollary \ref{cor01e}(a) and (b)
that if X eventually plays $\pi_X$ then against any Y play, the limit point set $\Omega$ is  contained in the
convex hull  $[W, (P,P),V, (R,R)] =  [W,(P,P),V,(R,R),(T,S)] \cap [W,(R,R),(S,T)]$.

Now assume that Y also eventually plays such a plan $\pi_Y$ with lines $\ell',\ell_1',\ell_2'= \ell_2$ and with points $V'$ and $W$.

Let $D$ consist of the set of points of $\S$ which are either on or below $)(P,P),(S,T)($ or on or below
$)(P,P),(T,S)($. Let $D_1 = [W,Switch(W),$ $(P,P)]$, the set of $s \in \S$ with $s_X, s_Y \leq P$. It is easy to check that
if $s^{N^*} \in D$ then $s^n \in D$ for all $n \geq N^*$ and that for some $N^{**} \geq N^*$, $s^{N^{**}} \in D_1$.
All subsequent outcomes are $dd$ and so $\lim s^N = (P,P)$. On the other hand, if $\{ s^N \}$ does not converge to
$(P,P)$ then
it is easy to check that some $s^{N^{**}}$ lies in $[(P,P),V,(R,R),Switch(V')] \setminus \{ (P,P) \}$. All subsequent
outcomes are $cc$ and so $\lim s^N = (R,R)$.  Thus, we always have either convergence to $(P,P)$ or to $(R,R)$.

$\Box$ \vspace{.5cm}

Return now to consider a simple Smale plan with separation line $\ell$. For such a plan the choices
on the line $\ell$ were, in general, left unspecified. Recall that if $\pi_X$ and $\pi_Y$ are simple Smale plans with separation
lines $\ell_X$ and $\ell_Y$ then, except for the extreme cases when both $\ell_X$ and $\ell_Y$ are the diagonal or
both the co-diagonal  (which requires $P \leq \frac{1}{2}(T+S)$), then if, eventually, X uses $\pi_X$ and Y uses
$\pi_Y$ the sequence $\{ s^N \}$ converges to the point of intersection  $\ell_X \cap Switch(\ell_Y)$ regardless of
earlier play and regardless of the choices on the separation lines.

To illustrate where the choices on the lines become important, let us consider the extreme cases.

Suppose $P < \frac{1}{2}(T+S)$ and $\pi$ is a simple Smale plan with separation line the co-diagonal,
$\ell  \ = \ )(S,T),(T,S)($. If both players use $\pi$ from time $N^*$ on
then if $s^{N^*}$ is above $\ell$ then the sequence remains above $\ell$ converging to $(T,S)$. Similarly, if
$s^{N^*}$ is below $\ell$ then the sequence remains below $\ell$ and converges to $(S,T)$.  If $s^{N^*} \in \ell$ then
the result depends on the choice of $\pi$ on $\ell$. If $\pi(s^{N^*}) = \pi(Switch(s^{N^*})) = 1$ then $s^{N^* +1}$ is
above $\ell$ with convergence to $(T,S)$. If $\pi(s^{N^*}) = \pi(Switch(s^{N^*})) = 0$,  then $s^{N^* +1}$ is
below $\ell$ with convergence to $(S,T)$. Suppose that $\pi(s) = 0$ if $s \in \ell$ with $s$  above the diagonal
and $\pi(s) = 1$ if $s \in \ell$ on or below the diagonal, then we again have alternating $(S,T)$ and $(T,S)$ motion
with limit $ \frac{1}{2}(T+S,T+S)$ unless the sequence lands on the point
$\frac{1}{2}(T+S,T+S)$ in which case we have convergence to $(T,S)$.

Now let $\pi$ be the simple Smale plan with separation line $\ell$ the diagonal.  Suppose that
for $s = (Q,Q) \in \ell$,  $\pi(s) = 1$ if $Q \geq \frac{1}{2}(T+S)$ and $= 0$ if $Q < \frac{1}{2}(T+S)$.
I think of this as the Smale version of Tit-for Tat. Suppose both players use $\pi$ for $N \geq N^*$ and
$s^{N^*} \not\in \ell$. We obtain alternating motion towards $(T,S)$ and $(S,T)$ with limit point
$\frac{1}{2}(T+S,T+S)$ unless at some time $N \geq N^*$, $s^N = (Q,Q) \in \ell$.
If $Q \geq \frac{1}{2}(T+S)$ then we obtain outcomes $cc$ for all rounds after $N$ with convergence to $(R,R)$.
If $Q < \frac{1}{2}(T+S)$ then we obtain outcomes $dd$ for all rounds after $N$ with convergence to $(P,P)$.

Thus, in both extreme cases, the limit results depend upon the $\pi$ choices on the separation line.
 \vspace{1cm}

 \section{Separation Paths and the Folk Theorem}
 \vspace{.5cm}

This section is the result of some suggestions and questions raised by Christian Hilbe in response to an earlier version.

\begin{df}\label{defnash01a}  We call a pair of Smale plans $\pi_X, \pi_Y$ a \emph{Nash Equilibrium}  when the following
 hold:

 \begin{itemize}
 \item[(a)] There is a point  $s^* = (s^*_X, s^*_Y)$ such that if X eventually plays  $\pi_X$ and Y eventually plays
 $\pi_Y$ then the outcome sequence $\{ s^N \}$ converges to $s^*$. We then call $s^*$ the \emph{payoff to the pair}.

 \item[(b)] If   Y eventually plays $\pi_Y$ then
  any limit point $V$ of an outcome sequence satisfies $V_X \leq s^*_X$, regardless of the play of X.\\
  If   X eventually plays $\pi_X$ then
  any limit point $V$ of an outcome sequence satisfies $V_Y \leq s^*_Y$, regardless of the play of Y.

  \end{itemize}
  \end{df}

  That is, neither player can obtain an improved payoff by  unilaterally changing plans.

  The so-called \emph{Folk Theorem of Iterated Play} in this context should say that for any $s^* \in \S$ with
  $s^*_X, s^*_Y \geq P$, there exists a Nash Equilibrium with payoff $s^*$.

  With  $R \geq s^*_X, s^*_Y \geq P$ this is easy to obtain.  If $\pi_X$ and $\pi_Y$ are the equalizer simple Smale plans with
  horizontal separation lines $\ell_X$ given by $s_Y= s^*_Y$ and $\ell_Y$ given by $s_Y = s^*_X$ then the pair
  $\pi_X, \pi_Y$ is a Nash equilibrium. In fact, as long as Y plays $\pi_Y$, the limit set $\Omega$ lies in the vertical
  line $Switch(\ell_Y)$ and X obtains $s^*_X$ regardless of his play
  Consequently, X has no incentive to use $\pi_X$. Similarly, for Y. So we would like to strengthen condition
  (b) to the analogue of the Protection Condition in this context.

 \begin{df}\label{defnash01b}  We call a pair of Smale plans $\pi_X, \pi_Y$ a \emph{Strong Nash Equilibrium}  when the following
 hold:

 \begin{itemize}
 \item[(a)] There is a point  $s^* = (s^*_X, s^*_Y)$ such that if X eventually plays  $\pi_X$ and Y eventually plays
 $\pi_Y$, then the outcome sequence $\{ s^N \}$ converges to $s^*$.

 \item[(b')] If  Y eventually plays $\pi_Y$, then the point $V = s^*$ is the only limit point $V$ of an outcome sequence with $V_X \geq s^*_X$, regardless of the play of X.\\
 If  X eventually plays $\pi_X$, then the point $V = s^*$ is the only limit point $V$ of an outcome sequence with $V_Y \geq s^*_Y$, regardless of the play of Y.

  \end{itemize}
  \end{df}

  Also, we would like to obtain Nash equilibrium results when $s^*_X$ or $s^*_Y$ is greater than $R$.

  To deal with the cases when $P > \frac{1}{2}(T + S)$, we again let $\bar P = \min(P,\frac{1}{2}(T + S))$.

  \begin{df}\label{defnash02} If $s \in \S$ then the \emph{upper triangle with vertex $s$}, denoted $T(s)$,
   is the triangle $[s,(S,T),(R,R)]$.
  The \emph{lower quadrangle with vertex $s$}, denoted $Q(s)$, is the convex set
  $[s,(T,S),(P,P),(\bar P, \bar P)]$.

  A non-empty subset $C \subset \S$ is \emph{upper  full} if $s \in C$ implies $T(s) \subset C$,
 and \emph{lower  full} if $s \in C$ implies $Q(s) \subset C$.
  \end{df}

  For $C \subset \S$ we let $C^{\circ}$ denote the $\S$ interior. Observe that if
  $s$ is on the segment $[(S,T),(R,R)]$ then $T(s)$ is the segment with empty interior. If $P \leq \frac{1}{2}(T + S)$
  so that $\bar P = P$, then $Q(s)$ is the triangle $[s,(T,S),(P,P)]$.

  If $C \subset \S$ is  nonempty, then the distance from $s$ to $C$ is $d(s,C) = \inf \{ ||s - \hat s|| : \hat s \in C \}$.
  For $\e \geq 0$, let $C_{\e} = \{ s \in \S : d(s,C) \leq \e \}$. Observe that $C_0$ is the closure of $C$.

  \begin{lem}\label{nash03} Let $C$ be a nonempty subset of $\S$ and let $\e \geq 0$.

  (a) If $C$ is upper full, then it is connected and contains $[(S,T),(R,R)]$. If $C$ is lower full,
  then it is connected and contains $[(T,S),(P,P),(\bar P, \bar P)]$.

  (b)  If $C$ is closed and is upper (or lower) full then $C_{\e}$ is upper full (resp. lower full).
  In particular, the closure of $C$ is upper full (resp. lower full).

  (c)  The union $C^+$ (or $C^-$) of the upper triangles (resp. the lower quadrangles) with vertices in $C$
   is an upper full (resp. lower full) subset
  which is closed if $C$ is.

  (d) If $C$ is convex then it is upper full if $(S,T),(R,R) \in C$ and it is
  lower full if  $(T,S),(P,P),(\bar P, \bar P) \in C$.
  \end{lem}

  {\bfseries Proof:} We prove the results for the upper case as the lower is similar.

  (a) If $C$ is upper full then it is a union of upper triangles, all of which are connected
  and all of which contain $[(S,T),(R,R)]$. Such a union  is connected.

  (b) Let $\e_1 > \e$. If $r \in \S$ then the segment $[s,r]$ consists of the points $ts + (1-t)r$ for all $t \in [0,1]$.
  The upper triangle $T(s)$ is the union of the segments $[s,r]$ with $r \in [(S,T),(R,R)]$.
  If $s_1 \in \S$ then $|| (ts + (1-t)r) - (ts_1 + (1-t)r) || = t || s - s_1 || \leq || s - s_1 ||$.
  Thus, if $s_1$ is $\e_1$ close to $s$, then every point of $T(s_1)$  is $\e_1$ close to a point
  of  $T(s)$  and vice-versa. It follows that if $C$ is upper full then $C_{\e}$ is upper full.
  The lower quadrangle $Q(s)$ is the union of the segments $[s,r]$ with $r \in [(T,S),(P,P),(\bar P, \bar P)]$.

  (c) If $s_1 \in T(s)$, then the $T(s_1) \subset T(s)$  by convexity of $T(s)$.
  Hence, any union of upper triangles is upper full.

  The set $C^+$
  is the image of the continuous map on $C \times [(S,T),(R,R)] \times [0,1]$ by
  $(s,r,t) \mapsto ts + (1-t)r$. If $C$ is closed, then the image is compact and so is closed. The set $C^-$
  is the image of the analogous continuous map on $C \times [(T,S),(P,P),(\bar P, \bar P)] \times [0,1]$

  (d) Obvious.

$\Box$ \vspace{.5cm}

Because upper and lower fullness are generalizations of convexity, the following is an extension of
Corollary  \ref{cor01d}.

\begin{prop}\label{propnash04} Assume that eventually X plays  a Smale plan $\pi$, Y uses an arbitrary plan and that
the initial play is arbitrary.  Let $\Omega$ be the limit point set of an associated sequence of outcomes. Let $C \subset \S$
be a nonempty closed set.

(a) If $C$ is lower full and  $ \S \setminus C$ is  contained in the defection zone $\pi^{-1}(0)$, then $\Omega \subset C$.

(b) If $C$ is upper full  and  $ \S \setminus C$  is  contained in  the cooperation zone $\pi^{-1}(1)$, then $\Omega \subset C$.
\end{prop}

{\bfseries Proof:} (a) Assume X plays $\pi_X$ after $N^*$.
For $\e > 0$ let $N_{\e} > N^* $ be an $\e$ step-time. If for all $N \geq N_{\e}$,
$s^N \in \S \setminus C$, then after $N_{\e}$ the outcomes are just as though X played All D. Hence, by
Proposition \ref{prop01c} (b), $\Omega \subset [(T,S),(P,P)] \subset C$.  So we may assume that at some time
$N_1 > N_{\e}$, $s^{N_1} \in C$. We show by induction that for all $N \geq N_1, \ s^N \in C_{\e}$.

If $s^N \in C$, then since $N > N_{\e}$, $|| s^{N+1} - s^N || < \e$, and so $s^{N+1} \in C_{\e}$.

If $s^N \in C_{\e} \setminus C \subset \pi^{-1}(0)$, then $S^{N+1} \in [(T,S),(P,P)]$ and so
$s^{N+1}$ is in $Q(s^N)$ and so is contained in the lower full set $C_{\e}$.

Since, eventually, $s^N$ is in the closed set $ C_{\e}$, it follows that $\Omega \subset C_{\e}$. As $\e$ was
arbitrary, $\Omega \subset \bigcap_{\e > 0} \ C_{\e} = C$.

The proof of (b) is completely analogous.

$\Box$ \vspace{.5cm}

\begin{df}\label{defnash05} A \emph{separation path} is a closed, connected subset $C \subset \S$ which meets
the segments $[(\bar P,\bar P),(S,T)]$ and $[(R,R),(T,S)]$ and which satisfies
\begin{itemize}
\item[(*)] If $s \in C$ then $T(s)^{\circ} \cup Q(s)^{\circ}$ is disjoint from $C$.
\end{itemize}
\end{df}
 \vspace{.5cm}

 Note that if $s \in [(S,T),(R,R)]$ then $T(s)$  has empty interior and so is disjoint from $C$.

 Recall that when $P > \frac{1}{2}(T + S)$ the line $)(S,T),(P,P)($ intersects $((R,R),(T,S))$ at a point
  we label $\bar W$. We let $\bar W$ denote the point $(T,S)$ when $P \leq \frac{1}{2}(T + S)$.
\vspace{.5cm}

 \begin{theo}\label{theonash06} Assume that  $C \subset \S$ is a separation path with $s \in C$. Let
 $proj: C \to [S,T]$ be the restriction to $C$ of the first coordinate projection, i.e. $proj(s) = s_X$.

\begin{equation}\label{eqnash1}
C^+ \cup C^- \ = \ \S \quad \text{and} \quad C^+ \cap C^- = C. \hspace{2cm}
\end{equation}
\begin{equation}\label{eqnash2}
\S \setminus C^- \ = \ C^+ \setminus C \ = \ (C^+)^{\circ} \quad \text{and}
\quad \S \setminus C^+ \ = \ C^- \setminus C \ = \ (C^-)^{\circ}.
 \end{equation}
\begin{equation}\label{eqnash3}
T(s)^{\circ} \cap C^- \ = \ \emptyset \quad \text{and} \quad Q(s)^{\circ} \cap C^+ \ = \ \emptyset.
 \end{equation}
\begin{equation}\label{eqnash4}
Q(\bar W)^{\circ} \cap C \ = \ \emptyset. \hspace{2cm}
 \end{equation}

 The map $proj$ is injective mapping $C$ onto an interval $[a,b]$ with $S \leq a \leq \bar P$ and $R \leq b \leq T$.
 \end{theo}

 {\bfseries Proof:}  Since $C$ meets $[(S,T),(\bar P,\bar P)]$ and $[(R,R),(T,S)]$, $C$
 meets $[V_-,V_+]$ if $V_-, V_+ = (S,T),(\bar P,\bar P)$ or
 if $V_-, V_+ = (R,R),(T,S)$. If $V_- \in ((\bar P,\bar P),(T,S))$ and $V_+ \in ((S,T),(R,R))$ then the
 segment $[V_-,V_+]$ separates the edges $[(S,T),(\bar P,\bar P)]$ and $[(R,R),(T,S)]$ and so meets the connected set $C$.

 Let $s_0 \in C \cap [V_-,V_+]$. The half-open  segment $(s_0,V_+]$ is contained in $T(s_0)^{\circ}$ and so
$(s_0,V_+] \subset C^+ \setminus C$ by Condition (*). Similarly, $[V_-,s_0) \subset C^- \setminus C$.
On the other hand, if $s_1 \in Q(s) \cap [V_-,V_+]$ then $[V_-,s_1) \subset Q(s_1)^{\circ} \subset  Q(s)^{\circ}$.
 So Condition (*) implies $s_0 \not\in [V_-,s_1)$. Thus,
$Q(s) \cap [V_-,V_+] \subset [V_-,s_0]$. Similarly, $T(s) \cap [V_-,V_+] \subset [s_0,V_+]$. Since $s$ was an arbitrary
point of $C$, including the possibility $s = s_0$, it follows
that $C^- \cap [V_-,V_+] = [V_-,s_0]$ and $C^+ \cap [V_-,V_+] = [s_0,V_+]$.
Thus, we have
$$[V_-,V_+]\cap (\S \setminus C^-) = [V_-,V_+]\cap ( C^+ \setminus C) = [V_-,V_+]\cap (C^+)^{\circ} = (s_0,V^+].$$
with analogous equations for $[V_-,s_0)$. Because every point of $\S$ lies on some such interval $[V_-,V_+]$,
equations  (\ref{eqnash2}) follow and clearly imply  (\ref{eqnash1}) . Since $T(s)^{\circ} \subset (C^+){\circ}$ and
$Q(s)^{\circ} \subset (C^-){\circ}$, (\ref{eqnash2}) implies (\ref{eqnash3}), as well.

If $P \leq \frac{1}{2}(T + S)$ then $\bar W = (T,S)$ and $Q(\bar W) = [(\bar P,\bar P),(T,S)]$ has empty interior.
$P > \frac{1}{2}(T + S)$ then for all $s \in Q(\bar W)^{\circ}$, $(P,P) \in T(s)^{\circ}$. Since $(P,P) \in C^-$,
(\ref{eqnash3}) implies that $s \not\in C$, proving (\ref{eqnash4}).

The image $proj(C)$ is a compact connected subset of $[S,T] \subset \R$. Hence, it is a closed interval $[a,b]$.
Since $C$ meets $[(S,T),(\bar P,\bar P)]$ and $[(R,R),(T,S)]$, $proj(C)$ meets the closed intervals $[S,\bar P]$ and $[R,T]$ in $\R$.

For $t \in \R$ with $S \leq t \leq T$ let $\ell_t$ be the vertical line $ \{ s_X = t \}$. For $t \in proj(C)$ with
 $S < t < R$  the
point $V_+' = \ell_t \ \cap \ )(S,T),(R,R)($ lies in $((S,T),(R,R))$. Let $s_t$ be the
point of $C \cap \ell_t$ with the smallest $s_Y$ coordinate. Since $(s_t,V_+'] \subset T(s_t)^{\circ}$ it follows
from Condition (*) that $C \cap \ell_t$ is the singleton $\{ s_t \}$. Similarly, for $t \in proj(C)$ with
$\bar P < t < T$,
$\ell_t$ meets $((\bar P,\bar P),(T,S))$ at a point $V_-'$. Let $s^t$ be the
point of $C \cap \ell_t$ with the largest $s_Y$ coordinate. $[V_-',s^t) \subset Q(s^t)^{\circ}$ and so
$C \cap \ell_t$ is the singleton $\{ s^t \}$. Finally, $(S,T)$ is the only point of $\S$ with $s_X = S$ and
$(T,S)$ is the only point of $\S$ with $s_X = T$. It follows that each vertical line meets $C$ in at most one point.
Thus, $proj$ is injective.

$\Box$ \vspace{.5cm}

If $\ell$ is a separation line, then $(S,T)$ and  $(R,R)$ are on or above $\ell$ and so for
$s \in \ell \cap \S$, $T(s) \subset H^+$ and $T(s)^{\circ} \subset H^+ \setminus \ell$.
Also, $(P,P)$ and $(T,S)$ are on or below $\ell$. Since $(R,R)$ is on or above $\ell$ it follows that
$(\bar P, \bar P)$ is also on or below $\ell$. Thus, $Q(s) \subset H^-$ and $Q(s)^{\circ} \subset H^- \setminus \ell$.
It follows that $\ell \cap \S$ is a separation path. Call $\ell$ a \emph{strict separation line} when
it is a separation line which does not contain $(S,T), (R,R), (T,S)$ or $(P,P)$. That is,  for a strict separation line $\ell$ the points $(S,T)$ and $(R,R)$ are strictly above $\ell$ and the points $(T,S)$ and $(P,P)$ are
strictly below $\ell$. Furthermore, $(\bar P, \bar P)$ is strictly below $\ell$ as well, because if
$ (\bar P, \bar P) \not= (P,P)$ then $P > \frac{1}{2}(T+S)$ and the only separation line with $ (\bar P, \bar P)$
on or above it is the diagonal which contains $(P,P)$. Thus, with $L$ an affine map associated with $\ell$,
we have $L(S,T), L(R,R) > 0$ and $L(T,S), L(P,P), L(\bar P, \bar P) < 0$. It follows that:

\begin{equation}\label{eqnash5}
s \in \ell \cap \S \quad \Longrightarrow \quad  (T(s) \cup Q(s)) \cap \ell = \{ s \}.
 \end{equation}

 By analogy, if $C$ meets
the segments $[(S,T),(\bar P,\bar P)]$ and $[(R,R),(T,S)]$ and satisfies the strengthening of Condition (*)
\begin{equation}\label{eqnash5a}
s \in C\quad \Longrightarrow \quad  (T(s) \cup Q(s)) \cap C = \{ s \},
 \end{equation}
 then we will call $C$ a \emph{strict separation path}.

Recall that a function $\g : [a,b] \to \R$ is called piecewise $C^1$ when $\g$ is continuous and there is a finite
sequence $a = a_0 <  a_1 < \dots < a_n = b$ so that $\g$ is continuously differentiable on each subinterval
$[a_{i-1},a_i]$ for $i = 1, \dots,n$. Then at each point $(t,\g(t))$ with $t \in [a,b)$ there is a tangent line from
the right and at each point with $t \in (a,b]$ there is a tangent line from the left. Except at the points
with $t = a_i$ for $i = 0, \dots, n$ these two agree are both are the true tangent line at the point.
\vspace{.5cm}

\begin{theo}\label{theonash06a} Assume that  $C$ is a subset of $\S$.
\begin{itemize}
\item[(a)] If $C$ is a separation path then there exists a continuous function $\g : [a,b] \to [S,T]$
with $S \leq a \leq \bar P, R \leq b \leq T$ such that $C$ is the graph of $\g$, i.e. $C = \{ (t,\g(t)) : t \in [a,b] \}$.

\item[(b)] If $C$ is a separation path which is the graph of a piece-wise $C^1$ function $\g$, then
every tangent line (i.e. every line tangent to a point from the left or the right) is a separation line.

\item[(c)] Assume that $C$ is the graph of a piece-wise $C^1$ function with $(a,\g(a)) \in [(\bar P,\bar P), (S,T)]$ and
$(b,\g(b)) \in [(R,R),(T,S)]$. If every tangent line is a strict separation line, then $C$ is a strict separation path.

\item[(d)] If $C_1$ and $C_2$ are separation paths, then $C_1 \cap Switch(C_2)$ is non-empty. If, in addition,
one of them is a strict separation path, then $C_1 \cap Switch(C_2)$ is a singleton set.
\end{itemize}
\end{theo}

{\bfseries Proof:} (a) As it is a continuous bijection on a compact set, $proj : C \to [a,b]$ is a homeomorphism.
If $\g$ is the composition of $proj^{-1}$ with the projection to the $s_Y$ coordinate then $\g$ is a continuous
map and $proj^{-1}$ is given by $t \mapsto (t,\g(t))$.

For (b) and (c) we use $s(t)$ for $(t,\g(t))$ with $t \in [a,b]$.

(b) If $t \in [a,b)$ then Condition (*) implies that for every $t_1 \in (t,b]$ the secant line $)s(t),s(t_1)($
passes through $[(R,R),(T,S)]$. Furthermore, if $t \leq P$, then $(P,P)$ is on or
below the secant line. Passing to the limit, the same is true for every tangent line from the right.
Similarly, for $t \in (a,b]$ every tangent line from the left
passes through $[(S,T),(\bar P, \bar P)]$ and if $t \geq P$ then
$(P,P)$ is on or below the tangent line from the left.

 Thus, at the points of $[a,b] \setminus \{ a_0,\dots,a_n \}$ the tangent line is a separation line.
The remaining tangent lines from the left and right are limits of true tangent lines and so are separation lines as well.

(c) The graph is a closed, connected set which meets $[(S,T),(\bar P,\bar P)]$ and $[(R,R),(T,S)]$ by hypothesis.
If $s \in Q(\bar W)^{\circ}$ then no separation line passes through $s$ and so the graph is disjoint from $Q(\bar W)^{\circ}$.

Let $a \leq t_0 < b$. Then (\ref{eqnash5}) for the tangent line from the right implies that, for sufficiently
small $h > 0$, $s(t_0+h) \not\in T(s(t_0)) \cup Q((s_0))$. If $s(t) \in T(s(t_0)) \cup Q((s_0))$ for some
$t > t_0$ then we can let $t^*$ be the first entrance time, i.e.
$t^* = \inf \{ t > t_0 : s(t) \in T(s(t_0)) \cup Q((s_0)) \}$ and let $\ell$ be the tangent line from
the left at $s(t^*)$. Thus, $s(t^*)$ is in the $\S$ topological
boundary of $T(s(t_0)) \cup Q((s_0))$. The portion of the boundary with $s_X > s_X(t_0)$ consists of
$(s(t_0),(R,R)]$ and either $(s(t_0),(T,S)]$ or if $(P,P)$ lies above this segment then
$(s(t_0),(P,P)] \cup ((P,P),(T,S)]$ (since $s(t_0) \not\in Q(\bar W)^{\circ}$). If $s(t^*) \in  (s(t_0),(R,R)]$
then because $s(t_1)$ is below
$(s(t_0),(R,R)]$ for $t_0 < t_1 < t^*$ the point $(R,R)$ lies below the secant line $)s(t_1),s(t^*)( $.
It follows that, in the limit, $(R,R)$ is on or below $\ell$. Similarly, if $s(t^*) \in (s(t_0),(T,S)]$
or $s(t^*) \in ((P,P),(T,S)]$ then $(T,S)$ is on or above $\ell$. Finally, if $s(t^*) \in (s(t_0),(P,P)]$
then $(P,P)$ is on or above $\ell$.

Let $a < t_0 \leq b$. Then, for sufficiently
small $h > 0$, $s(t_0-h) \not\in T(s(t_0)) \cup Q((s_0))$.  If $s(t) \in T(s(t_0)) \cup Q((s_0))$ for some
$t < t_0$ then we can let $t^*$ be the first entrance time moving left, i.e.
$t^* = \sup \{ t < t_0 : s(t) \in T(s(t_0)) \cup Q((s_0)) \}$ and let $\ell$ be the tangent line from
the right at $s(t^*)$. The portion of the boundary with $s_X < s_X(t_0)$ consists of
$(s(t_0),(S,T)]$ and either $(s(t_0),(\bar P,\bar P)]$ or if $(P,P)$ lies above this segment then
$(s(t_0),(P,P)] \cup ((P,P),(\bar P,\bar P)]$(since $s(t_0) \not\in Q(\bar W)^{\circ}$).
As before $s(t^*) \in (s(t_0),(S,T)]$ implies
$(S,T)$ is on or below $\ell$, and $s(t^*) \in (s(t_0),(\bar P,\bar P)]$ or $s(t^*) \in ((P,P),(\bar P,\bar P)]$
implies $(\bar P, \bar P)$ is on or above $\ell$. Finally,  if $s(t^*) \in (s(t_0),(P,P)]$ then $(P,P)$ is
on or above $\ell$.

In none of these cases can $\ell$ be a strict separation line.  This proves (\ref{eqnash5a})
 for any $t = t_1$. This in turn implies implies $C$ is a separation path.

(d) The connected set $Switch(C_2)$ meets $[(S,T),(R,R)] \subset C_1^+$ and
$[(\bar P, \bar P),(T,S)] \subset C_1^-$. It follows from
(\ref{eqnash1}) for $C_1$ that $Switch(C_2)$ meets $C_1$.

For any $s \in \S$, $(T(s) \cup Q(s)) \cup Switch((T(s)^{\circ} \cup Q(s)^{\circ}) = \S$. Assume that
$C_1$ is a strict separation path and that $s_1 \in \S \setminus \{ s \}$. If $s_1 \in C_1$ then
$s_1 \not\in (T(s) \cup Q(s))$. If $s_1 \in Switch(C_2)$ then $s_1 \not\in Switch((T(s)^{\circ} \cup Q(s)^{\circ})$.
It follows that $s_1 \not\in C_1 \cap Switch(C_2)$.

$\Box$ \vspace{.5cm}

{\bfseries Remark:} It is possible to extend (c) somewhat.  It is easy to check that the set of separation paths
is itself closed in the space of closed subsets of $\S$ equipped with the Hausdorff topology. So if $\{\g_n \}$
is a sequence of piecewise $C^1$ functions which satisfy the conditions of (c) and the sequence of graphs converges
in the appropriate sense to the graph of $\g$ then the graph of $\g$ is a separation path.
\vspace{.5cm}

Using simple geometric arguments similar to those in (c) above,
we can describe what occurs in a non-strict separation path. Since we
will not need the results, we leave the proof to the interested reader.
\vspace{.5cm}

\begin{prop}\label{propnash06b} Assume that $s_1$ and $s_2$ are distinct points of a separation path $C$.
\begin{itemize}
\item[(i)] Let $V_-,V_+$ be one of the pairs $(S,T),(R,R)$ or $(P,P),(T,S)$ or $(\bar P, \bar P), (P,P)$.
If $C$ meets the open segment $(V_-,V_+)$, then it contains the closed segment $[V_-,V_+]$.

\item[(ii)] Let $V$ be one of the points $(S,T),(R,R),(T,S),(P,P)$ or $(\bar P, \bar P)$. If
$s_2 $ lies in the open segment $(s_1,V)$ then $C$ contains the closed segment $[s_1,s_2]$.
\end{itemize}
\end{prop}
%
%
%

$\Box$ \vspace{.5cm}

\begin{ex} A differential equations construction for strict separation paths. \end{ex}

For the example, we will restrict to the case $P < \frac{1}{2}(T+S)$and we will just sketch the argument,
leaving the details to the reader.

Define for $s \in \S$
\begin{equation}\label{eqnash5b}
\begin{split}
m^+(s) \ = \ \begin{cases} (R - s_Y)/(R - s_X) \quad \text{if} \ s_Y \geq s_X \\
(s_Y - P)/(s_X - P) \quad \text{if} \ s_Y \leq s_X. \end{cases}  \hspace{.5cm}\\
m^-(s) \ = \ \begin{cases} (T - s_Y)/(s_X - S) \quad \text{if} \ s_Y + s_X \geq T + S\\
(s_Y - S)/(T - s_X) \quad \text{if} \ s_Y + s_X \leq T + S. \end{cases} \\
\end{split}
\end{equation}
So $m^+(s) = 1 $ if $s_Y = s_X$ and $m^-(s) = 1$ if $s_Y + s_X = T + S$. If $s \in [(S,T),(R,R)]$ or
$s \in [(P,P),(T,S)]$ then $m^+(s) = -m^-(s)$. Otherwise, $m^+(s) > -m^-(s)$.
A line $\ell $ through $s$ with slope $m$ is a separation line iff
$m^+(s) \geq m \geq -m^-(s)$ and it is a strict separation line iff both inequalities are strict.

Now let $m$ be a  real-valued, differentiable function  on $\S$ with $m(s) = m^+(s) = -m^-(s)$ on
$[(S,T),(R,R)] \cup [(P,P),(T,S)]$ and with $m^+(s) > m(s)  > -m^-(s)$ otherwise.

Consider the differential equation $\frac{dy}{dx} = m(x,y)$ defined for $s = (x,y) \in \S$. Observe that
$[(S,T),(R,R)]$ and $[(P,P),(T,S)]$ are both solution curves for the differential equation. So by the uniqueness
theorem for ode's every other solution curve remains in $\S \setminus ([(S,T),(R,R)] \cup [(P,P),(T,S)])$
and so by   Theorem \ref{theonash06a} (c) each of the remaining solution curves is a strict separation path.

$\Box$ \vspace{.5cm}

Since the $\S$ intersection with a separation line is a separation path, the following extends Corollary \ref{cor01b}.
\vspace{.5cm}

\begin{cor}\label{cornash07} Assume that eventually X plays  a Smale plan $\pi$, Y uses an arbitrary plan and that
the initial play is arbitrary.  Let $\Omega$ be the limit point set of an associated sequence of outcomes.
If $C$ is a separation path such that
$C^+ \setminus C \subset \pi^{-1}(0)$ and $C^- \setminus C \subset \pi^{-1}(1)$, then $\Omega \subset C$. \end{cor}

{\bfseries Proof:} $C^+$ is upper full with $\S \setminus C^+ = C^- \setminus C$ contained in the cooperation zone and
$C^-$ is lower full with $\S \setminus C^- = C^+ \setminus C$ contained in the defection zone. It follows from
Proposition \ref{propnash04} that $\Omega \subset C^+ \cap C^- = C$.

$\Box$ \vspace{.5cm}

\begin{prop}\label{propnash08} If $s^* = (s^*_X,s^*_Y) \in \S$ with $P < s^*_Y < R$ then there is
a strict separation path $C$ such that $s^*$ is the unique point  $s \in C$ with $s_Y \geq s^*_Y$. \end{prop}

{\bfseries Proof:} First, assume $s^* \not\in [(S,T),(\bar P,\bar P)] \cup [(R,R),(T,S)]$.

Choose $V$ a point on the open segment $((S,T),(\bar P,\bar P))$ with
$P < V_Y < s^*_Y$ and so that  $)V,s^*($ intersects $((R,R),(T,S))$. Choose $W$ a point on the open segment $((R,R),(T,S))$
with $P < W_Y < s^*_Y$ and so that  $)s^*,W($ intersects $((S,T),(\bar P,\bar P))$. The lines $)V,s^*($ and $)s^*,W($
are strict separation lines.

Let $C = [V,s^*] \cup [s^*,W]$. While it is easy to check directly that $C$ satisfies (\ref{eqnash5a}), it follows from
Theorem \ref{theonash06a}(c) that $C$ is a strict separation path. Clearly, $s^*$ is the unique point of $C$ with maximum height.

If $s^* \in [(S,T),(\bar P,\bar P)]$ then choose $W$ as above so that $\ell = )s^*,W($ is a strict separation line.
 If $s^* \in [(R,R),(T,S)]$ then choose $V$ as above so that $\ell = )s^*,V($ is a
 strict separation line. In either case, $s^*$ is
 the point of maximum height on the segment $\ell \cap \S$.

$\Box$ \vspace{.5cm}

\vspace{.1cm}
\includegraphics{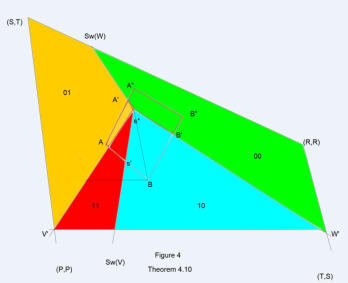}
\vspace{.1cm}

\begin{theo}\label{theonash09} If $s^* = (s^*_X,s^*_Y) \in \S$ with $P < s^*_X,  s^*_Y $ then there is
a pair of Smale plans $\pi_X, \pi_Y$ which is a strong Nash equilibrium with $s^*$ the payoff to the pair. \end{theo}

{\bfseries Proof:}  \textbf{Case 1} ($P < s^*_X,  s^*_Y < R$):  Apply Proposition \ref{propnash08} to choose $C_X$
a strict separation path with maximum point $s^*$
 and $C_Y$ be a strict separation path with maximum
point $Switch(s^*) = (s^*_Y,s^*_X)$. Let $\pi_X$ and $\pi_Y$ be Smale plans with
\begin{equation}\label{eqnash11a}
\begin{split}
C_X^+ \setminus C_X \subset (\pi_X)^{-1}(0), \quad C_X^- \setminus C_X \subset (\pi_X)^{-1}(1), \\
C_Y^+ \setminus C_Y \subset (\pi_Y)^{-1}(0), \quad C_Y^- \setminus C_Y \subset (\pi_Y)^{-1}(1).
\end{split}
\end{equation}
Corollary \ref{cornash07} and Theorem \ref{theonash06a} (d) imply that when $\pi_X$ plays
$\pi_Y$ the payoff is $s^* = C_X \cap Switch(C_Y)$.
If X uses an alternative plan, then
the limit set is contained in $Switch(C_Y)$, and similarly if Y varies against $\pi_X$.
 \vspace{.5cm}

\textbf{Case 2} ($s^*_X = s^*_Y = R$): This is the good plan case studied in Section 3. Use $\pi_X = \pi_Y$ a simple
Smale plan with the separation line through $(R,R)$ with slope $m$ satisfying $0 < m < 1$.
 \vspace{.5cm}

\textbf{Case 3} ($P < s^*_X < R \leq s^*_Y$ , or $P < s^*_Y < R \leq s^*_X$): Notice that if $s \in \S$ with
$s_Y > R$, then $s_X < R$. Furthermore,  for the case $P < s^*_Y < R \leq s^*_X$ it suffices to apply $Switch$ to
the other case.  So we will assume that $P < s^*_X < R \leq s^*_Y$. To begin with we will also assume that
$s^* \not\in [(R,R),(S,T)]$.

First we use Proposition \ref{propnash08},
or more precisely its proof, to choose $C_Y$  a separation path with maximum
point $Switch(s^*) = (s^*_Y,s^*_X)$.  That is, $C_Y = [Switch(V),Switch(s^*)] \cup [Switch(s^*),Switch(W)]$ with
$V \in ((\bar P,\bar P),(T,S)), W \in ((R,R),(S,T))$, with $P < V_X, W_X < s^*_X$ and with the other conditions as described above.
Let $\ell = )V,s^*($.
Let $\pi_Y = 0$ on $C_Y^+ \setminus C_Y$ and $ = 1$ on $C_Y^-\setminus C_Y$. So regardless of the X plan, the Switched
limit set $Switch(\Omega)$ is contained in $C_Y$ which has $Switch(s^*)$ as it unique point of maximum height.

Next choose $W' \in ((R,R),(T,S))$  and  $V' \in [(S,T),(\bar P,\bar P)]$
 with $P < W'_Y = V'_Y < R $.  So the horizontal
line $\ell' = )V',W'($ is a separation line. Let $C_X = [V',s^*] \cup [s^*,W']$. Let $\pi_X = 0$ above
$C_X$ and $= 1$ below $C_X$. Now despite the labeling, $C_X$ is not a separation path. The lower full set
$C_X^-$ is the set of points on or below each of the lines $)V',s^*( $ and $)s^*,W'($,
but the set of points on or above $C_X$ is not
upper full. For every point $r \in [V',s^*)$, $(r,s^*] \subset T(r)^*$.
On the other hand, Corollary \ref{cor01b} (c) applies to $\pi_X$ with $C \subset C_X^-$ equal to the triangle $ [V',s^*,W']$.
 So regardless of the Y plan, the
limit set $\Omega$ is contained in $C$ which has $s^*$ as its unique point of maximum height.

To complete the proof we must show that if, after some $N^*$, X uses $\pi_X$ and Y uses $\pi_Y$ then the solution
sequence converges to $s^*$, i.e. $\Omega = \{ s^* \}$. We know that $\Omega $ is contained in
$C \cap Switch(C_Y)$.  If $\bar s = \ell \cap \ell'$,  then $C \cap Switch(C_Y)$ is the segment $K = [s^*,\bar s]$.
Thus, given an arbitrary $\d > 0$ there is an $N^{\d} > N^*$ after which that the sequence  remains $\d$ close to $K$.

Let $\d > 0$ be arbitrary. Within the $\d$ ball  $V_{\d}(s^*)$ we will build a box which contains $s^*$ in its interior.
Choose $A \in (V',s^*) \cap V_{\d}(s^*) $
close enough to $s^*$ that the line $\ell_1 = )(S,T),A($ crosses $\ell$ within $V_{\d}(s^*)$,
i.e. $s' = \ell_1 \cap \ell \in V_{\d}(s^*)$. Choose
$B \in \ell_1 \cap V_{\d}(s^*)$ so that $s' \in (A,B) \subset V_{\d}(s^*)$. We may choose $A$ close enough that
the horizontal segment from $B$ to $(s^*,V')$ is also in $V_{\d}(s^*)$. Now let $\ell_2 = )(P,P),B($
and $\ell_3 = )(P,P),A($.  Because $(P,P)$ lies below $)V',s^*($ and above $\ell = )V,s^*($, it follows that
$s^*$ is below $\ell_3$ and above $\ell_2$. We can choose $A$ close enough to $s^*$ that $B' = \ell_2 \cap (s^*,W')$
and $A' = \ell_1 \cap (s^*,W)$ lie in $V_{\d}(s^*)$.  This is where we use $s^* \not\in [(R,R),(S,T)]$. Choose
$B'' \in \ell_2 \cap V_{\d}(s^*), A'' \in \ell_1 \cap V_{\d}(s^*)$ so that $B' \in (B,B'')$ and $A' \in (A,A'')$.
The region of $\S$ above both $)s^*,W'($ and $)s^*,W($ is convex, as is $V_{\d}(s^*)$. It follows that the
quadrilateral $Box = [A,B,B'',A'']$ lies in $V_{\d}(s^*)$ and contains $s^*$ in its interior.

The various lines
cut $\S$ into a -finite- number of regions and we choose $\e < \d$
to be smaller than the distances between any two
disjoint such regions and small enough that the open neighborhood $V_{\e}(K)$ is contained in the union of $Box$ and
$C^-$, the convex set of points of $\S$ in $C$ or below $\ell'$ and, finally, small enough that $V_{\e}(Box)$ is contained
in $V_{\d}(s^*)$.

 Let $N_{\e} > N^*$ be an $\e$ step-time so that for every $N \geq N_{\e}$
the length $|| s^{N+1} - s^N||$ is less than $\e$. Hence, no such small move crosses between disjoint regions.
Since $\Omega \subset K$ we can also choose $N_{\e}$ so that
for every $N \geq N_{\e}$ $s^N \in V_{\e}(K)$.

Suppose that for some  $N  \geq N_{\e}$, $s^N \not\in Box$. If $s^N$ lies to the left of $\ell$ then
$S^{N+1} = (R,R)$ and if $s^N$ lies to the right of $\ell$ then
$S^{N+1} = (S,T)$. For $s^N \in \ell$ either $S^{N+1} = (R,R)$ or $S^{N+1} = (S,T)$.
For any such $s$ the (negative) slope of $)(S,T),s($ is less than the slope of $)(R,R),s($.
Thus  from the left we move toward $(R,R)$ and after crossing we move toward $(S,T)$. The successive crossings
are higher on $K$ and thus our net motion is upward until we enter $Box$. From  $Box$ it is easy to see that exit
could only occur from  the triangle $[s',B,s^*]$ landing above the horizontal line through $B$ and $\e$ close to
$Box$ and so is still in $V_{\d}(s^*)$. If it lands to the left of $\ell$ then the motion toward $(R,R)$ is
closer to $Box$ and so $[s^N,s^{N+1}] \subset V_{\d}(s^*)$. If it lands to the right of $\ell$ then
the motion toward $(S,T)$ is
upward, remaining above the $B$ horizontal line and so $[s^N,s^{N+1}] \subset V_{\d}(s^*)$. Subsequent alternating
moves are above these initial ones and so remain in $ V_{\d}(s^*)$ until the sequence re-enters $Box$. It follows that
eventually the sequence lies in  $V_{\d}(s^*)$ and so the limit point set $\Omega$ is contained
in the closed ball of radius $\d$ about $s^*$. As $\d > 0$ was arbitrary, it follows that $\Omega = \{ s^* \}$
as required.

Finally, we adjust the proof to deal with the case $s^* \in [(S,T),(R,R)]$. If $s^* = (R,R)$ we are in Case 2.
Since $s^*_X > P$, we have $s^* \not= (S,T)$. We use $C_X$ as before. Now $C_Y$ is the separation line
$\ell = )s^*,V($ with $W$ undefined. The points $A, B, B'$ and the lines $\ell_1, \ell_2, \ell_3$ are chosen as before.
This time, $B'' = \ell_2 \cap )(S,T),(R,R)($ and $A'' = \ell_3 \cap )(S,T),(R,R)($. Now $s^*$ is in the
upper edge $(A'',B'')$ which is nonetheless in the $\S$ interior of $Box  = [A,B,B'',A'']$. With an easy adjustment of
the previous argument, one can again show that $\Omega$ is contained
in the closed ball of radius $\d$ about $s^*$ and so conclude that $\Omega = \{ s^* \}$.

$\Box$ \vspace{1cm}

\section{Competition Among Simple Smale Plans}
 \vspace{.5cm}

In this section we move beyond the classical question which motivated our original interest in good strategies. We
consider now the evolutionary dynamics among simple Smale plans.  We follow Hofbauer and Sigmund \cite{HS} Chapter 9 and
Akin \cite{A-90}.

The dynamics that we consider takes place in the context of a symmetric two-person game, but generalizing our initial
description, we merely assume that there is a  set  of strategies indexed by a finite set $\I$.
 When players X and Y use strategies
with index $i, j \in \I$,
respectively, then the payoff to player X is given by $A_{ij}$ and the payoff to Y is $A_{ji}$. Thus, the game is
described by the payoff matrix $\{ A_{ij} \}$. We imagine a population of players each using a particular strategy for each
encounter and let $\xi_i$ denote the ratio of the number of $i$ players to the total population.  The frequency
vector $\{ \xi_i \}$ lives in the unit simplex $\Delta \subset \R^{\I}$, i.e. the entries are nonnegative and sum to $1$.
The vertex $v(i)$ associated with $i \in \I$ corresponds to a population consisting entirely of $i$ players. Thus,
$\xi = v(i)$ exactly when $\xi_i = 1$.  We assume the
population is large so that we can regard $\xi$ as changing continuously in time.

Now we regard the payoff in units of \emph{fitness}.  That is, when an $i$ player meets a $j$ player in an interval of
time $dt$, the payoff $A_{ij}$ is an addition to the background reproductive rate $\rho$ of the members of the population. So
the $i$ player is replaced by $1 + (\rho + A_{ij})dt  \ i$ players. Averaging over the current population distribution,
the expected relative reproductive rate for the subpopulation of $i$ players is $\rho + A_{i \xi}$, where
\begin{equation}\label{32}
\begin{split}
A_{i \xi} \quad = \quad \Sigma_{j \in  \ \I} \ \xi_j A_{ij} \qquad \mbox{ and} \hspace{2cm}\\
A_{\xi \xi} \quad = \quad \Sigma_{i \in  \ \I} \ \xi_i A_{i \xi}  \quad = \quad \Sigma_{i,j  \ \in \ \I} \ \xi_i \xi_j A_{ij}.
\end{split}
\end{equation}

The resulting dynamical system on $\Delta$ is given by the \emph{Taylor-Jonker Game Dynamics Equations} introduced in
Taylor and Jonker \cite{TJ}.

\begin{equation}\label{33}
\frac{d \xi_i}{dt} \quad = \quad \xi_i (A_{i \xi} \ - \ A_{\xi \xi} ).
\end{equation}

This system is an
example of the \emph{replicator equation} systems studied in great detail in Hofbauer and Sigmund \cite{HS}.

We will need some
general game dynamic results for later application. Fix the game matrix $\{ A_{ij} \}$.

A subset $A $ of $\Delta$ is called \emph{invariant} if $\xi(0) \in A$ implies that the entire solution
path lies in $C$.  That is, $\xi(t) \in A$ for all $t \in \R$. An invariant point is an
\emph{equilibrium}.

Each nonempty subset $\J$ of $ \I$ determines the \emph{face} $\Delta_{\J}$ of the simplex consisting of those $\xi \in \Delta$
such that $\xi_i = 0 $ for all $i \not\in \J$.
Each face of the simplex is  invariant  because $\xi_i = 0$ implies
that $ \frac{d \xi_i}{dt} = 0$. In particular,
 for each $i \in \I$ the vertex $v(i)$, which represents fixation at the $i$ strategy, is an equilibrium.

 In general,
 $\xi$ is an equilibrium when, for all $i, j \in \I$, $\xi_i, \xi_j > 0$ imply $A_{i \xi} = A_{j \xi}$, or, equivalently,
 $A_{i \xi} = A_{\xi \xi}$ for all $i$ such that $\xi_i > 0$, i.e. for all $i$ in the
 \emph{support} of $\xi$.

 An important example of an invariant set is the \emph{omega limit point set of an orbit}. Given an initial point
  $\xi \in \Delta$ with
 associated solution path $\xi(t)$, it is defined by intersecting the closures of the tail values.
 \begin{equation}\label{omega}
 \omega(\xi) \quad = \quad \bigcap_{t > 0} \overline{ \{ \xi(s) : s \geq t \}}.
 \end{equation}
 By compactness this set is nonempty. A point is in $\omega(\xi)$ iff it is the limit of some
sequence $\{ \xi(t_n) \}$ with $\{ t_n \}$ tending to infinity.
The set  $\omega(\xi)$ consists of a single point $\xi^*$ iff $Lim_{t \to \infty} \xi(t) = \xi^*.$
In that case, $\{ \xi^* \}$ is an invariant point, i.e. an equilibrium.

Notice that this is the analogue for the solution path of the limit point set $\Omega$ of a payoff sequence, considered
in the previous sections.
\vspace{.5cm}

\begin{df}\label{def-ess} We call a strategy $i^*$ a
\emph{evolutionarily stable strategy} (hereafter, an ESS) when
\begin{equation}\label{ess}
A_{j i^*} \ < \ A_{i^* i^*} \qquad \mbox{for all} \ \ j \not= i^* \quad \mbox{in} \ \ \I. \hspace{2cm}
\end{equation}
\end{df}
 \vspace{.5cm}

 {\bfseries Remark:} We follow \cite{A-16} in labeling this condition ESS although it is stronger than the condition
 originally introduced in \cite{M-S}. In  \cite{HS}  precisely this condition is called a \emph{strict Nash equilibrium} and so  this language requires a bit of justification.

 In the pure game theory context we regard a distribution $\xi$ over $I$ as a mixed-strategy rather than a population
 distribution of pure strategists. Then a pair $\xi_1, \xi_2$ is a \emph{Nash equilibrium} when each is a best reply against the
 other. That is, for all distributions $\eta$,
 \begin{equation}\label{nash1}
A_{\eta \xi_2} \ \leq \ A_{\xi_1 \xi_2} \quad \mbox{and} \quad A_{\eta \xi_1} \ \leq \ A_{\xi_2 \xi_1}.
\end{equation}
Notice that from this we see that
 \begin{align}\label{nash2}
  \begin{split}
  (\xi_1)_j > 0 \quad &\Longrightarrow \quad A_{j \xi_2} \ = \ A_{\xi_1 \xi_2} \qquad \mbox{and}\\
    (\xi_2)_j > 0 \quad &\Longrightarrow \quad A_{j \xi_1} \ = \ A_{\xi_2 \xi_1}.
\end{split}
\end{align}
That is, all of the pure strategies active in $\xi_1$ are best replies to $\xi_2$ and vice-versa.
In the context of Smale strategies this is the concept used in the previous section.

Following \cite{A-90} we call the pair  a \emph{regular Nash equilibrium}  when, in addition to (\ref{nash1}),
 \begin{align}\label{nash3}
 \begin{split}
(\xi_1)_j = 0 \quad &\Longrightarrow \quad A_{j \xi_2} \ < \ A_{\xi_1 \xi_2}  \qquad \mbox{and} \\
(\xi_1)_j = 0 \quad &\Longrightarrow \quad A_{j \xi_2} \ < \ A_{\xi_1 \xi_2}.
\end{split}
\end{align}
That is, the pure strategies active in $\xi_1$ are all of the pure strategies which give the best reply to $\xi_2$ and
vice-versa.

Returning to the dynamic context, a distribution $\xi$ is called a (regular) Nash equilibrium when the
pair $\xi, \xi$ is a (regular) Nash equilibrium in the above sense. From (\ref{nash2}) we see that a Nash equilibrium is
an equilibrium for the Taylor-Jonker equations as defined above. When $\xi$ is the vertex $v(i^*)$ then
it is a regular Nash equilibrium  exactly when it is a strict Nash equilibrium as defined on page 62 of \cite{HS},
or, equivalently, (\ref{ess}) holds.
 \vspace{.5cm}

\begin{prop}\label{prop-ess} If $i^*$ is an ESS then the vertex $v(i^*)$
 is an attractor, i.e. a locally stable equilibrium, for the system (\ref{33}).  In fact, there exists $\e > 0$ such that
\begin{equation}\label{34a}
1 \ > \ \xi_{i^*} \ \geq \ 1 - \e \qquad \Longrightarrow \qquad  \frac{d \xi_{i^*}}{dt} \ > \ 0.
\end{equation}
Thus, near the equilibrium $v(i^*)$, which is characterized by $\xi_{i^*} = 1, \ \xi_{i^*}(t)$
increases monotonically, converging to $1$ and the
alternative strategies are eliminated from the population in the limit.
\end{prop}

{\bfseries Proof:}  When $i^*$ is an ESS, $ A_{i^* i^*} > A_{j i^*}$ for all $j \not= i^*$.
It then follows for $\e > 0$ sufficiently small that
$  \xi_{i^*}   \geq  1 - \e$ implies $A_{i^* \xi} > A_{j \xi}$ for all $j \not= i^*$. If also $1 > \xi_{i^*}, $ then
$A_{i^* \xi} > A_{\xi \xi}$. So (\ref{33}) implies (\ref{34a}).

$\Box$ \vspace{.5cm}

\begin{df}\label{def-dom} For $\J$ a nonempty subset of $ \I $ we say a strategy $i$
\emph{weakly dominates} a strategy $j$  in $\J$ when $i, j \in \J$ and
\begin{equation}\label{dom}
A_{j k} \ \leq \ A_{i k} \qquad \mbox{for all} \ \ k  \in  \ \J, \hspace{2cm}
\end{equation}

with strict inequality for $k = i$  or $k = j$.   If the inequalities are strict for all $k$ then
we say that $i$ \emph{dominates} $j$ in $\J$.

We say that $i \in \J$ weakly dominates a sequence $\{ j_1, ..., j_n \}$ in $\J$ when there
exists $1 \leq m \leq n$ such that  $i$ weakly dominates $j_p$ in $\J$ for $ p = 1,\dots, m$ and
for $p = m+1, ..., n$, $i$ dominates $j_p$ in $\J \setminus \{ j_1,..., j_{p - 1} \}$.

When $\J$ equals all of $\I$ we will omit the phrase ``in $\J$".
\end{df}
\vspace{.5cm}

\begin{prop}\label{prop-dom} For $i \in \I$, let $\xi(t)$ be a solution path with $\xi_i(0)  > 0$

(a) If $i$ weakly dominates $j$ then $Lim_{t \to \infty} \ \xi_j(t) \  = \  0.$

(b) If $i$ weakly dominates the sequence $\{ j_1,...,j_n \}$ then for
$j = j_1, ..., j_n,$  $ Lim_{t \to \infty} \ \xi_j(t) \quad = \quad 0. $
\end{prop}

{\bfseries Proof:} (a): The face $\{ \xi : \xi_j = 0 \}$ is invariant. So if $\xi_j(0) = 0$ then $\xi_j(t) = 0$ for all $t$ and
so it is $0$ in the limit.  Thus, we may assume $\xi_j(0) > 0$.

For $i, j \in \I$, define the open set $Q_{ij}$ and on it the real valued function $H_{ij}$ by
\begin{equation}\label{Lij}
\begin{split}
Q_{ij} \quad = \quad \{ \xi \in \Delta : \xi_i, \xi_j > 0 \}  \\
H_{ij}(\xi) \quad = \quad \ln(\xi_i) - \ln(\xi_j).
\end{split}
\end{equation}

Let $h_0 = H_{ij}(\xi(0))$.

Observe that on $Q_{ij}$
\begin{equation}\label{Lijeq}
 dH_{ij}/dt \ = \ A_{i \xi} - A_{j \xi} \ = \
\Sigma_{k \in  \ \I} \xi_k (A_{ik} - A_{jk}) > 0.
\end{equation}

Hence, $ H_{ij}(\xi(t))$ is a strictly increasing function of $t$
on the open invariant set $Q_{ij}$. Thus, as a $t$ tends to infinity $H_{ij}(\xi(t))$ approaches
$h_{\infty} = sup \{ H_{ij}(\xi(t)) : t \geq 0 \}$ with
$h_0 < h_{\infty} \leq + \infty$.

We must prove  that $\xi_j = 0$ on the omega limit set.
Assume instead that $\xi^* \in \omega(\xi(0))$ with $\xi^*_j > 0$.  If $\xi^*_i$ were $0$ then $H_{ij}(\xi(t))$ would
not be bounded below on $\{ \xi(t) : t \geq 0 \}$.  Hence, $\xi^*$  lies in
$Q_{ij}$ with $h_{\infty} = H_{ij}(\xi^*) < \infty$. So on the invariant set
$\omega(\xi(0)) \cap Q_{ij}$, which contains $\xi^*$, and so is nonempty, $H_{ij}$ would be constantly $h_{\infty} < \infty$. Since
this set is invariant, $ dH_{ij}/dt $ would equal zero. This contradicts
(\ref{Lijeq}) which implies that the derivative is positive on $\omega(\xi(0)) \cap Q_{ij}$.

The proof of (b) is a  variation of the proof of (a). We refer to  \cite{A-16} Proposition 4.6.  An obvious  adjustment
of the initial step in the inductive proof of (b) there yields the proof here.

 $\Box$ \vspace{.5cm}

 \begin{cor}\label{cor-dom} Assume $I = \{ i^*,j_1, \dots, j_n \}$ and $i^* \in I$
 weakly dominates the sequence $\{ j_1,...,j_n \}$.
If $\xi_{i^*}(0) > 0$ then
 $\lim_{t \to \infty}  \ \xi_{i^*}(t) = 1$.
 \end{cor}

 {\bfseries Proof:} By Proposition \ref{prop-dom} $\xi_{j_p}(t) \to 0$ for all $p = 1,\dots, n$ and so
 $\xi_{i^*}(t) = 1 - \Sigma_{p = 1}^n \xi_{j_p}(t) \to 1$.

 $\Box$ \vspace{.5cm}

In \cite{A-16}, see also \cite{A-13}, we examined competition
among certain special Markov plans called \emph{Zero-Determinant Plans}.
It was proved that good Markov plans among them are attractors when competing against plans which are not agreeable.
In addition, global stability was proved when the class of competitors was further restricted.  Here we will
similarly consider competition among simple Smale plans. Let $I = \{ i^*,j_1, \dots, j_n \}$ index a list
of simple Smale plans with $\pi_i$ associated with separation line $\ell_i$ for $i \in I$. Except for the extreme
cases the intersection  $\ell_i \cap Switch(\ell_j)$ is a single point. For $\pi$ the plan with $\ell$ the diagonal
we will assume that the plan is weakly agreeable and that the initial play is $c$. So if X and Y both play
$\pi$ the payoff is $(R,R)$. If $P \leq \frac{1}{2}(T + S)$ then the co-diagonal is a separation line
and we will adopt the convention that if both players use co-diagonal plans then the payoff is
$\frac{1}{2}(T + S,T + S)$. If X plays $\pi_i$ and Y plays $\pi_j$ we will let $(A_{ij},A_{ji})$ be
the coordinates of the payoff point with the above conventions in the extreme cases. We will then use (\ref{33})
to represent the dynamics of the competition with $\xi_i$ the fraction of the $\pi_i$ players in the population.

Notice that if Y plays an equalizer plan $\pi_j$ with $\ell_j$ horizontal then $A_{ij} = A_{i'j}$ for
any plans $\pi_i$ and $\pi_{i'}$ for X. In particular, if all of the plans are equalizer plans then
$A_{i\xi} = A_{i'\xi }$ for all $i, i' \in I$ and so $A_{i \xi} = A_{\xi \xi}$ for all $i \in I$ and for any
population state $\xi$. Thus, the dynamics is trivial and every state $\xi$ is an equilibrium.

Now we consider the case when $\ell_{i^*}$ is a protection line.  That is, $\ell_{i^*}$ is a line through $(R,R)$ with
slope $m$ satisfying $0 < m \leq 1$. Notice that $m = 1$ is the diagonal line case. If $m < 1$ and
$\pi_{i^*}(R,R) = 1$ then $\pi_{i^*}$ is
a good simple Smale plan.
\vspace{.5cm}

\begin{theo}\label{theo-ess} If $\ell_{i*}$ is a protection line and $(R,R) \not\in \ell_j$ for any $j \in I \setminus \{ i^* \}$
then $i^*$ is an ESS and so fixation at $i^*$ is an attractor. \end{theo}

{\bfseries Proof:} $A_{i^* i^*} = R$. If Y plays $\pi_{i^*}$ and X plays $\pi_j$ for $j \in I \setminus \{ i^* \}$
then the payoff point is not $(R,R)$ and so $A_{j i^*}$ and $A_{i^* j}$ are both less than $R$ because $\ell_{i^*}$
is a protection line. This implies (\ref{ess})
and so the result follows from Proposition \ref{prop-ess}.

 $\Box$ \vspace{.5cm}

Since $\xi_{i^*}(0) = 0$ implies $\xi_{i^*}(t) = 0$ for all $t$ the best stability result we can hope for
is that every solution with $\xi_{i^*}(0) > 0$ converges to fixation at $i^*$. We will call this
\emph{global stability}.

As an illustration we describe a very special case.
\vspace{.5cm}

\begin{theo}\label{theo-equalizer} Assume that $\ell_{i^*}$ is a protection line. If for every
$j \in I \setminus \{ i^* \}$, $\ell_j$ is a horizontal line $\{ s_Y = C_j \}$ with $P \leq C_j < R$,
then for every $j \in I \setminus \{ i^* \}$, $i^*$ weakly dominates $j$ and so the system
exhibits global stability. \end{theo}

{\bfseries Proof:} Because $\ell_{i^*}$ is a protection line, we have, as in Theorem \ref{theo-ess},
$A_{j i^*} < A_{i^* i^*} = R$. For any $k \in I \setminus \{ i^* \}$, $C_k = A_{jk} = A_{i^*k}$ for all
$j \in I$ and weak domination, (\ref{dom}), follows.  A fortiori, $i^*$ weakly dominates the sequence
$\{j_1,\dots,j_n \}$ and the result follows from Corollary \ref{cor-dom}.

 $\Box$ \vspace{.5cm}

We will show that we achieve global stability if $\pi_{i^*}$ is good,
$(R,R) \not\in \ell_j$ for any $j \in I \setminus \{ i^* \}$  and, in addition, all the lines $\ell_j$ have positive
slope. This requires a bit of geometry.
\vspace{.5cm}

\begin{lem}\label{lem-glo} Assume that $\ell_{i^*}$ is a protection line,
$(R,R) \not\in \ell_k$ for any $k \in I \setminus \{ i^* \}$
 and that  $\ell_k$ has non-negative slope for every $k \in I$. If for some
$\bar i \in I $ the segment $\ell_{\bar i} \cap \S  $
lies below $\ell_{i^*}$ then $i^*$ weakly dominates
$\bar i$. \end{lem}

{\bfseries Proof:} As usual $A_{\bar i i^*}  < A_{i^* i^*} = R$.  For any $k \in I$ the line $Switch(\ell_k)$ is
either vertical or has positive slope. Let $\bar V, V^*$ be the intersection points of $Switch(\ell_k) \cap
\ell_{\bar i}$ and  $Switch(\ell_k) \cap \ell_{i^*}$, respectively. If $Switch(\ell_k)$ is vertical then the X coordinates
of $\bar V$ and $V^*$ are equal.  If $Switch(\ell_k)$ has positive slope then $V^*$ is above and to the right of
$\bar V$ and so has a larger X coordinate. Thus, $A_{\bar i k} \leq A_{i^* k}$, proving weak domination.

 $\Box$ \vspace{.5cm}

Now we assume that the slope of $\ell_{i^*} < 1$, i.e. $\ell_{i^*}$ is not the diagonal. Let $V$ be the point of
intersection  $\ell_{i^*} \ \cap \ ) (S,T),(\bar P,\bar P) ( $, where $\bar P = \min(P, \frac{1}{2}(T + S))$.
Thus, $\ell_{i^*} \cap \S \setminus {(R,R)} = [ V, (R,R) )$ and the entire half-open segment lies above the diagonal.
Now let $\ell_j$ be a separation line  which does not contain $(R,R)$.
So it contains a point $A \in ((R,R),(T,S)]$.
Let $B$ be the intersection point  $\ell_j \ \cap \ ) (S,T),(\bar P,\bar P) ( $. If $B$ lies below
$\ell_{i^*}$ then the entire segment $\ell_j \cap \S$ lies below  $\ell_{i^*}$. Otherwise, $B \in [V,(S,T)]$ and this
is the situation we wish to examine.

Since $A$ is below $\ell_{i^*}$ and $B$ is on or above $\ell_{i^*}$ it follows that $\ell_j$ intersects
$\ell_{i^*}$ at a point $V^j$ of $\S$ with $V^j_X$  its X coordinate. Notice that the
portion of $\ell_{j}$ to the right of $V^j$ lies below $\ell_{i^*}$. In any case,
$Switch(\ell_j)$ intersects $\ell_{i^*}$ at a point $W^j$ with X coordinate  $W^j_X$.
\vspace{.5cm}

\begin{lem}\label{geom} If $\ell_j$ has non-negative slope then $V^j_X < W^j_X$. \end{lem}

{\bfseries Proof:} The lines $\ell_j$ and $Switch(\ell_j)$ meet the diagonal at a common point
$(Q,Q) = \ell_j \ \cap \ Switch(\ell_j)$. To the right of $\{ s_X = Q \}$ the line $\ell_j$ lies below the
diagonal, because $A$ is below the diagonal. On the other hand, all of
$\ell_{i^*} \ \cap \ \S \setminus (R,R)$ lies above the diagonal. Hence, $V^j_X < Q$.
Similarly, $Switch(\ell_j)$ intersects $\ell_{i^*}$ above the diagonal. Since $Switch(\ell_j)$ is either
vertical or has positive slope, it follows that $Q \leq W^j_X$.


 $\Box$ \vspace{.5cm}

From this we obtain the main result of this section.
\vspace{.5cm}

\begin{theo}\label{theo-glo} $\{ \pi_i : i \in I \}$ be a finite indexed collection of simple Smale plans
with $\ell_i$ the separation line for $\pi_i$. Assume that
for some $i^* \in I, \ \ell_{i*}$ is a line through $(R,R)$ with slope strictly between $0$ and $1$
and $(R,R) \not\in \ell_j$ for any $j \in I \setminus \{ i^* \}$.  If  $\ell_i$ has non-negative slope for all $i \in I$
and, in addition,  $\ell_i \cap \S$ lies below $\ell_{i^*}$ for those $i \in I$ with $\ell_i$ horizontal, then
fixation at $i^*$ is a globally stable equilibrium. That is, if $\xi_{i^*}(0) > 0$ then
$\lim_{t \to \infty} \ \xi_{i^*}(t) = 1$. \end{theo}

{\bfseries Proof:} We choose a numbering of the $n$ elements of $I \setminus \{i^*\}$ by letting $j_1,...,j_m$
with $0 \leq m \leq n$ so that $\ell_j \cap \S$ lies below $\ell_{i^*}$ if and only if $j = j_p$ for
some $p \leq m$. If no such exist then $m = 0$ and the set is empty.

For the remaining $\ell_j$'s the slope is positive and the numbers $V^j_X$ and $W^j_X$ are
defined as above. Number them so that $V^{j_p}_X \leq V^{j_{p+1}}_X$ for $m < p < n$.

By Corollary \ref{cor-dom} it suffices to show that $i^* $ weakly dominates the sequence
$\{j_1, \dots, j_n \}$.

To begin with $i^*$ weakly dominates each $ j_p$ for $p \leq m$ by Lemma \ref{lem-glo}.

We must show that if $m < p \leq n$ then $i^*$ dominates $j_p$ in\\  $\{ i^*, j_p, j_{p+1},\dots, j_n \}$.

As before, $A_{j_p i^*} < A_{i^* i^*}$. Now let $k \in  \{  j_p, j_{p+1},\dots, j_n \}$.

Because of
the chosen numbering and Lemma \ref{geom} we have $V^{j_p}_X \leq V^k_X < W^k_X$. That is,
intersection point $W^k$ of $Switch(\ell_k) \ \cap \ \ell_{i^*}$ lies to the right of $V^{j_p}$.
The slope of  $Switch(\ell_k)$ is greater than $1$ and the slope of $\ell_{i^*}$ is less than $1$.
Hence, $Switch(\ell_k)$ is above $\ell_{i^*}$ to the right of $W^k$ and below $\ell_{i^*}$ to the left.
It follows that $Switch(\ell_k)$ intersects the vertical line  $\{ s_X = V^{j_p}_X \}$ below $V^{j_p}$
and so below the line $\ell_{j_p}$, because $V^{j_p}$ lies on $\ell_{j_p}$.
Again  $Switch(\ell_k)$ has slope greater than $1$ and $\ell_{j_p}$ has slope less than one.
So  $Switch(\ell_k)$ intersects $\ell_{j_p}$ to the right of
this vertical line. Right of this $V^{j_p}$ vertical line, $\ell_{j_p}$ lies below $\ell_{i^*}$. As in Lemma \ref{lem-glo}
the intersection point $Switch(\ell_k) \ \cap \ \ell_{j_p}$ lies
below and to the right of $W^{k} = Switch(\ell_k) \ \cap \ \ell_{i^*}$
That is, $A_{j_p k} < A_{i^* k}$. Thus, $i^*$ dominates $j_p$ in  $\{ i^*, j_p, j_{p+1},\dots, j_n \}$, as required.

 $\Box$ \vspace{1cm}

\section{Variations}
\vspace{.5cm}

Following Smale we consider alternative weighting schemes.

Let $\{ w_1, w_2, \dots \}$ be an infinite sequence of positive numbers. Let $W_N = \Sigma_{k=1}^N \ w_k$ and
$\Delta_N = \Sigma_{k = 1}^{N} |w_{k+1} - w_k|$. Consider the conditions:
\begin{itemize}
\item (Weight Condition 1) \quad $\lim_{N \to \infty} \  \frac{w_N}{W_N} = 0$.
\item (Weight Condition 2) \quad $\lim_{N \to \infty} \  W_N \ = \ \infty$.
\item (Weight Condition 3) \quad $\lim_{N \to \infty}  \frac{\Delta_N}{W_N} = 0$.
\end{itemize}

\begin{lem} Condition (3) implies Conditions (1) and (2).  If the sequence $\{ w_n \}$ is
monotonically non-increasing or non-decreasing, then Conditions (1) and (2) imply Condition (3).
\end{lem}

{\bfseries Proof:} $\Delta_N \geq |w_{N+1} - w_1|$. So (3) implies $\frac{ |w_{N+1} - w_1|}{W_N} \to 0$.
If $ |w_{n+1} - w_1| \leq \frac{1}{2} w_1$ infinitely often then $w_{n+1} > \frac{1}{2} w_1$ infinitely often
and so the increasing sequence $\{ W_N \}$ is unbounded, implying (2). Otherwise, eventually
$|w_{n+1} - w_1| > \frac{1}{2} w_1$ and so $\frac{w_1}{2 W_N} \to 0$ which implies (2).
Then $w_{n+1} \leq w_1 + |w_{n+1} - w_1|$ and so $\frac{w_{N+1}}{W_{N+1}} < \frac{w_{N+1}}{W_{N}} \to 0$ which is
(1).

Condition (1) implies $\frac{W_N}{W_{N+1}} = 1 - \frac{w_{N+1}}{W_{N+1}} \to 1$. Hence, (1) implies
$\frac{w_{N+1}}{W_N} \to 0$. Condition (2) implies $\frac{w_1}{W_N} \to 0$.

If the sequence $\{ w_n \}$ is monotone then the sum defining $\Delta_N$ telescopes to yield
$\Delta_N = | w_{N+1} - w_1 |$. So in that case (1) and (2) imply (3).

 $\Box$ \vspace{.5cm}

If the sequence is non-increasing then (1) certainly holds and (2) says that the
sequence does not decrease so fast that the associated series converges.
If the sequence is non-decreasing then (2) certainly holds and (1) says  that the
sequence does not increase too fast.  For example, if $w_{N+1} \geq \e W_N$ then
$\frac{w_{N+1}}{W_{N+1}} \geq \frac{\e}{1 + \e}$.

The initial averaging procedure that we used had $w_n = 1$ for all $n$ and so $W_N = N$.

Since the averaging procedure uses ratios we may multiply by a positive constant and so assume $w_1 = 1$ and
hence $W_N \geq 1$ for all $N$.

Now assume that $\{ w_n \}$ is a positive sequence with $w_1 = 1$ and Conditions (1) and (2) hold.

We replace our previous averaging of the payoff sequence in (\ref{4}) to define
\begin{equation}\label{av4}
s^N \quad = \quad \frac{1}{W_N} \ \Sigma_{k=1}^{N} \ w_k S^k.
\end{equation}

We obtain the analogues of (\ref{4a}) and (\ref{4b}).
\begin{equation}\label{av4a}
s^{N+1} \quad = \quad \frac{W_N}{W_{N+1}} s^N \ + \ \frac{w_{N+1}}{W_{N+1}} S^{N+1}.
\end{equation}
and so
\begin{equation}\label{av4b}
s^{N+1} - s^N \quad = \quad  \frac{w_{N+1}}{W_{N+1}} (S^{N+1} - s^N).
\end{equation}

By Condition (1) (\ref{av4}) implies that $|| s^{N+1} - s^N || \to 0$ and so the limit point set is connected as before.
However, the crucial fact is (\ref{av4a}) which says that $s^{N+1}$ is on the segment $[S^{N+1},s^N]$ with the
weight on $s^N$ approaching $1$ as $N \to \infty$. Consequently, all of the linear estimates for Smale plans
go through as before. The only change is that the numerical estimates $M N^*/N$ are replaced by
$M W_{N^*}/W_N$ which tends to $0$ as $N \to \infty$ by Condition (2). In particular, when two
non-extreme, simple Smale plans compete
we obtain convergence to the intersection point regardless of the averaging procedure.

It is the similar result for Markov plans that requires Condition (3).

Suppose that $\MM$ is the Markov matrix when X plays $\pp$ and Y plays $\qq$.  Let $\vv^1$ be the initial
distribution and $\vv^{n+1} = \vv^n \MM$, the distribution after round $n+1$.  Define
\begin{equation}\label{Mar4}
\bar \vv^N \quad = \quad \frac{1}{W_N} \ \Sigma_{k=1}^{N} \ w_k \vv^k.
\end{equation}
It follows that
\begin{equation}\label{Mar4a}
\bar \vv^N \MM \  = \  \frac{1}{W_N} \ \Sigma_{k=1}^{N} \ w_k \vv^{k+1} =   \frac{1}{W_N} \ \Sigma_{k=2}^{N+1} \ w_{k-1} \vv^{k}.
\end{equation}

Since the length of a distribution is at most $1$ we have that
\begin{equation}\label{Mar4b}
|| \bar \vv^N - \bar \vv^N \MM || \leq \frac{w_1 + w_{N+1} + \Delta_N}{W_N}.
\end{equation}

From Condition (3) it follows that any limit point of the sequence $\{\bar \vv^N \}$ is a stationary distribution.
In particular if there is a unique terminal set and so a unique stationary distribution $\vv$ then $\{\bar \vv^N \}$
converges to $\vv$.  If $J$ is one of several terminal sets then with probability $p_J$, depending only on then initial
distribution, $\vv^1$, the sequence of outcomes enters $J$. The conditional distributions assuming entrance into $J$
then converge to the unique stationary distribution on $J$.

In contrast with all this, there is another sort of natural averaging which does not work. Suppose we use
\begin{equation}\label{av4new}
s^N \quad = \quad \frac{1}{W_N} \ \Sigma_{k=1}^{N} \ w_k S^{N+1-k}.
\end{equation}

With Condition (3) one can still show that $|| s^{N+1} - s^N || \to 0$, but this time $s^{N+1}$ is not on the
segment $[S^{N+1}, s^N]$ except when all the $w_n$ are equal, in which case the two sorts of averaging agree (This is
the original $w_n = 1$ for all $n$ case).  So the results from the first section will not carry over.
\vspace{.5cm}

The other variation to consider is a asymmetric version of the Prisoner's Dilemma with payoffs given by
\renewcommand{\arraystretch}{1.5}
\begin{equation}\label{asym1}
\begin{array}{|c||c|c|}\hline
X\backslash Y & \quad c \quad & \quad  d \quad \\ \hline \hline
c & \quad (R_X,R_Y) \quad & \quad (S_X,T_Y) \quad  \\ \hline
d & \quad (T_X ,S_Y)\quad  & \quad (P_X,P_Y) \quad \\ \hline
\end{array}
\end{equation}
and with inequalities for X and for Y analogous to those of (\ref{3}).

This is a real issue because in the classic version of the Prisoner's Dilemma the payoffs are not in units of
dollars, time reduced from a prison sentence or population fitness, but in terms of utility and there is no reason
that the two players would have the same Von Neumann-Morgenstern utility functions.

At first glance, there is no problem.  In \cite{A-15} the good Markov strategies are characterized for the asymmetric case.
In \cite{Sm} Smale points out that the theory will work the same way for the asymmetric case.  Now one must describe separate
Smale strategies for Y, rather than using $\pi \circ Switch$, but as he indicates the mathematics is essentially the same.

There is, however, an underlying philosophical problem. In \cite{A-15} the inequalities for a good plan for X use the payoffs
for Y, which, in theory, X does not know. In the Markov case, this is not too bad because only a rough estimate is needed
to ensure that the strategy is good.

In the Smale case, the running averages use the payoffs to both players. Perhaps the best way to proceed would be to
begin again and operate, not in the two dimensional convex set generated by the payoff pairs but in the three dimensional
simplex of outcomes.  That is, let
\begin{equation}\label{asym2}
\begin{split}
\ee_{cc} = (1, 0, 0, 0), \qquad \ee_{cd} = (0, 1, 0, 0), \\
\ee_{dc} = (0, 0, 1, 0), \qquad \ee_{dd} = (0, 0, 0, 1).
\end{split}
\end{equation}
The convex hull $\S'$ with these vertices is the simplex of distributions on the four outcomes. The data we use from the
sequence of outcomes $\{ o_1, \dots, o_N \}$ is the frequency of past outcomes:
\begin{equation}\label{asym3}
s^N = \frac{1}{N} \Sigma_{k = 1}^N \ee_{o_k}.
\end{equation}
so that, analogous with (\ref{4a})
\begin{equation}\label{asym5}
s^{N+1} = \frac{1}{N+1} o_{N+1} \ + \ \frac{N}{N+1} s^N.
\end{equation}
A plan for X is then a map $\pi : \S' \to [0,1]$ with $\pi(s)$ the probability of cooperating in response to position $s$.
So a pure strategy plan, of the sort Smale uses would be a map $\pi : \S' \to \{ 0, 1 \}$.

Linear results analogous to those of Section 2 can then be carried over.
Nonetheless, determining what is a good plan would still require some estimate
of the opponent's payoffs. This is a task for another day.

\vspace{1cm}

\bibliographystyle{amsplain}

\end{document}